\documentclass[12pt,a4paper]{article}
\usepackage{epsfig,latexsym,amsfonts,amssymb,amsmath,amscd,
graphics,theorem,epic}
\theoremstyle{plain}
\newtheorem{lemma}{Lemma}
\newtheorem{proposition}[lemma]{Proposition}
\newtheorem{remark}[lemma]{Remark}

\newtheorem{theorem}{Theorem}

\newtheorem{corollary}[lemma]{Corollary}
{\theorembodyfont{\rmfamily} 
\font\ncsc=cmcsc10  \font\ntt=cmtt12
\begin{document}
\baselineskip=15pt
\newcommand{\pperp}{\hbox{$\perp\hskip-6pt\perp$}}
\newcommand{\ssim}{\hbox{$\hskip-2pt\sim$}}
\newcommand{\N}{{\mathbb N}}\newcommand{\Delp}{{\Pi}}
\newcommand{\A}{{\mathbb A}}
\newcommand{\Z}{{\mathbb Z}}
\newcommand{\R}{{\mathbb R}}
\newcommand{\C}{{\mathbb C}}
\newcommand{\Q}{{\mathbb Q}}\newcommand{\T}{{\mathbb T}}\newcommand{\K}{{\mathbb K}}
\newcommand{\PP}{{\mathbb P}}
\newcommand{\st}{{*}}
\newcommand{\mnote}{\marginpar}\newcommand{\ev}{{\operatorname{ev}}}
\newcommand{\Id}{{\operatorname{Id}}}\newcommand{\irr}{{\operatorname{irr}}}\newcommand{\nod}{{\operatorname{nod}}}
\newcommand{\oeps}{{\overline\eps}}\newcommand{\Area}{{\operatorname{Area}}}\newcommand{\End}{{\operatorname{End}}}
\newcommand{\oDel}{{\widetilde\Del}}
\newcommand{\real}{{\operatorname{Re}}}
\newcommand{\conv}{{\operatorname{conv}}}
\newcommand{\Span}{{\operatorname{Span}}}
\newcommand{\Ker}{{\operatorname{Ker}}}
\newcommand{\Fix}{{\operatorname{Fix}}}
\newcommand{\sign}{{\operatorname{sign}}}
\newcommand{\Log}{{\operatorname{Log}}}
\newcommand{\oi}{{\overline i}}
\newcommand{\oj}{{\overline j}}
\newcommand{\ob}{{\overline b}}
\newcommand{\os}{{\overline s}}
\newcommand{\oa}{{\overline a}}
\newcommand{\oy}{{\overline y}}
\newcommand{\ow}{{\overline w}}
\newcommand{\ou}{{\overline u}}
\newcommand{\ot}{{\overline t}}
\newcommand{\oz}{{\overline z}}
\newcommand{\newi}{i}
\newcommand{\newj}{j}
\newcommand{\newm}{m}
\newcommand{\newl}{{\ell}}
\newcommand{\bw}{{\boldsymbol w}}\newcommand{\bi}{{\boldsymbol i}}
\newcommand{\bx}{{\boldsymbol p}}\newcommand{\bp}{{\boldsymbol p}}
\newcommand{\bpp}{{\boldsymbol P}}
\newcommand{\by}{{\boldsymbol q}}
\newcommand{\bz}{{\boldsymbol z}}
\newcommand{\eps}{{\varepsilon}}
\newcommand{\proofend}{\hfill$\Box$\bigskip}
\newcommand{\Int}{{\operatorname{Int}}}
\newcommand{\pr}{{\operatorname{pr}}}
\newcommand{\grad}{{\operatorname{grad}}}
\newcommand{\rk}{{\operatorname{rk}}}
\newcommand{\im}{{\operatorname{Im}}}
\newcommand{\sk}{{\operatorname{sk}}}
\newcommand{\const}{{\operatorname{const}}}
\newcommand{\Sing}{{\operatorname{Sing}}}
\newcommand{\conj}{{\operatorname{Conj}}}
\newcommand{\Pic}{{\operatorname{Pic}}}
\newcommand{\Crit}{{\operatorname{Crit}}}
\newcommand{\Ch}{{\operatorname{Ch}}}
\newcommand{\discr}{{\operatorname{discr}}}
\newcommand{\Tor}{{\operatorname{Tor}}}
\newcommand{\Conj}{{\operatorname{Conj}}}
\newcommand{\val}{{\operatorname{val}}}
\newcommand{\Val}{{\operatorname{Val}}}
\newcommand{\res}{{\operatorname{res}}}
\newcommand{\add}{{\operatorname{add}}}
\newcommand{\tmu}{{\C\mu}}
\newcommand{\ov}{{\overline v}}\newcommand{\on}{{\overline n}}
\newcommand{\ox}{{\overline{x}}}
\newcommand{\tet}{{\theta}}
\newcommand{\Del}{{\Delta}}
\newcommand{\bet}{{\beta}}
\newcommand{\kap}{{\kappa}}
\newcommand{\del}{{\delta}}
\newcommand{\sig}{{\sigma}}
\newcommand{\alp}{{\alpha}}
\newcommand{\Sig}{{\Sigma}}
\newcommand{\Gam}{{\Gamma}}
\newcommand{\gam}{{\gamma}}
\newcommand{\Lam}{{\Lambda}}
\newcommand{\lam}{{\lambda}}
\newcommand{\SC}{{SC}}
\newcommand{\MC}{{MC}}
\newcommand{\nek}{{,...,}}
\newcommand{\cim}{{c_{\mbox{\rm im}}}}
\newcommand{\mathto}{\mathop{\to}}
\newcommand{\op}{{\overline p}}

\newcommand{\w}{{\omega}}

\title{A Caporaso-Harris type formula for Welschinger invariants
of real toric Del Pezzo surfaces}
\author{Ilia Itenberg \and Viatcheslav Kharlamov \and Eugenii Shustin}

\date{}
\maketitle

\begin{abstract}
We define a series of relative tropical Welschinger-type
invariants of real toric surfaces. In the Del Pezzo case, these
invariants can be seen as real tropical analogs of relative
Gromov-Witten invariants, and are subject to a recursive formula.
As application we obtain new formulas for Welschinger invariants
of real toric Del Pezzo surfaces.


\end{abstract}

\section{Introduction}\label{intro}

Welschinger invariants of real rational symplectic four-manifolds
\cite{W,W1} represent one of the most interesting and intriguing
objects in real enumerative geometry. In the case of a real
unnodal ({\it i.e.}, not containing any rational $(-n)$-curve,
$n\ge 2$) Del Pezzo surface $\Sig$ the Welschinger invariants
count, with appropriate weights $\pm 1$, the real rational curves
which belong to an ample linear system~$|D|$ and pass through a
given generic conjugation-invariant set of $c_1(\Sig)\cdot D-1$
points in~$\Sig$. In this paper we consider only the invariants
corresponding to sets of real points.

Our goal is to provide recursive formulas
which calculate the Welschinger invariants of toric Del Pezzo
surfaces equipped with the tautological real structure.
The formulas we obtain are similar
to those proved by L.~Caporaso and J.~Harris~\cite{CH}
for relative
Gromov-Witten invariants of $\PP^2$.

We use the technique of tropical geometry and follow ideas of
A.~Gathmann and H.~Markwig~\cite{GM1,GM2} who suggested a tropical
version of the Caporaso-Harris formula and its tropical proof. We
suitably adapt the tropical count to the real setting and
introduce tropical, multi-component and irreducible, Welschinger
numbers for relative constraints and arbitrary genus. We check
their invariance (Theorem~\ref{cht1} in Section~\ref{chsec2}) and
prove that these tropical invariants satisfy Caporaso-Harris type
formulas (Theorem \ref{formula} for the multi-component invariants
and Theorem \ref{cht3} for irreducible invariants,
Section~\ref{chsec3}; a reformulation with generating functions is
presented in Section~\ref{chsec5}). In the case when the set of
relative constraints is empty these invariants coincide with the
genuine Welschinger invariants. As a by-product, we establish some
monotonicity of the Welschinger invariants and give a new proof of
their positivity (Corollaries~\ref{chc3} and~\ref{chc1},
Section~\ref{chsec4}).

The paper is organized as follows. In Section 2 we remind
definitions and basic facts concerning Welschinger invariants and
plane tropical curves. Section 3 contains the definition of
tropical relative Welschinger numbers, the statements of the main
results and few corollaries. We prove the invariance of tropical
relative Welschinger numbers in Section 4 and the recursive
formulas in Section 5. Section 6 is devoted to concluding remarks.
Our main results are stated in terms of embedded plane tropical curves,
while the proofs go essentially through the parameterized
incarnation of these curves, and starting from Section 4 we put a special
attention to be maximally possible consistent with existing in the
literature notions and statements concerning these two different
categories. In particular, a part of Section 4 is devoted to various types of
genericity conditions and their comparison.

To conclude this short introduction, we would like to emphasize a
certain, challenging in our opinion, difference between the real
and complex cases. Namely, Gathmann-Markwig's count of tropical
curves in the tropical version of Caporaso-Harris formula is in a
strict correspondence (in the sense of modified Mikhalkin's
correspondence theorem \cite{Mi}, see \cite{GM1,GM2}) with the
count of complex algebraic curves in the original Caporaso-Harris
formula. In particular, the invariance of the terms in
Caporaso-Harris formula explains (and implies) the invariance of
the terms in the tropical version of Caporaso-Harris formula
proposed by Gathmann and Markwig. On the contrary, we do not know
how to lift up invariantly the terms entering the formulas
suggested in the present paper (except those which lead to the
genuine Welschinger invariants). One of the difficulties is that
such a lift, if it exists, can not be formulated in purely
topological terms, see section \ref{remarks}.

\medskip

{\bf Acknowledgements}. A considerable part of this work was done
during our visits to the Max-Planck-Institut f\"{u}r Mathematik,
Bonn. We thank MPIM for the hospitality and excellent working
conditions.

The authors were partially supported by a grant from the Ministry
of Science and Technology, Israel, and Minist\`{e}re des Affaires
Etrang\`{e}res, France. The first two authors were partially
funded by the ANR-05-0053-01 grant of Agence Nationale de la
Recherche and a grant of Universit\'{e} Louis Pasteur, Strasbourg.
The first and the third author enjoyed a support from the
Hermann-Minkowski-Minerva Center for Geometry at the Tel Aviv
University. The third author acknowledges a support from the grant
no. 465/04 from the Israeli Science Foundation.

\section{Preliminaries}\label{invariants}

\subsection{Welschinger invariants
of Del Pezzo surfaces}\label{Welsch} We remind here the definition
of Welschinger invariants~\cite{W, W1} restricting ourselves to a
particular situation. Let $\Sig$ be a real {\it unnodal} ({\it
i.e.}, not containing any rational $(-n)$-curve, $n\ge 2$) Del Pezzo
surface with a connected real part $\R\Sigma$, and let
$D\subset\Sigma$ be a real ample divisor. Consider a generic
set~$\boldsymbol\omega$ of $c_1(\Sig)\cdot D-1$ real points of
$\Sig$. The set $R(D, \boldsymbol\omega)$ of real $C\in|D|$ passing
through the points of~$\boldsymbol\omega$ is finite, and all these
curves are nodal and irreducible. (In fact, the listed properties of
$R(D, \boldsymbol\omega)$ can be taken here as a definition of the
term `generic'.) Due to the Welschinger theorem~\cite{W, W1} (and
the genericity of the complex structure on $\Sig$), the number
$$W
(\Sig,D, \boldsymbol\omega)=\sum_{C\in R(D,
\boldsymbol\omega)}(-1)^{s(C)}\ ,$$ where $s(C)$ is the number of
{\it solitary nodes} of $C$ ({\it i.e.}, real points, where a local
equation of the curve can be written over $\R$ in the form
$x^2+y^2=0$), does not depend on the choice of a generic
set~$\boldsymbol\omega$. We denote this Welschinger invariant by
$W(\Sig,D)$.

\subsection{Divisors on toric Del Pezzo
surfaces}\label{ample}
There are five
toric unnodal Del Pezzo surfaces: the projective plane
$\PP^2$, the product $\PP^1 \times \PP^1$
of projective lines, and $\PP^2$ with
$k$ blown up generic points, $k=1,2$, or 3; the latter three surfaces
are denoted by $\PP^2_k$.
Let $E_1, \ldots ,E_k$ be the
exceptional divisors of $\PP^2_k \to \PP^2$ and $L\subset \PP^2_k$
the pull back of a generic straight line.

\begin{figure}
\begin{center}
\epsfxsize 125mm \epsfbox{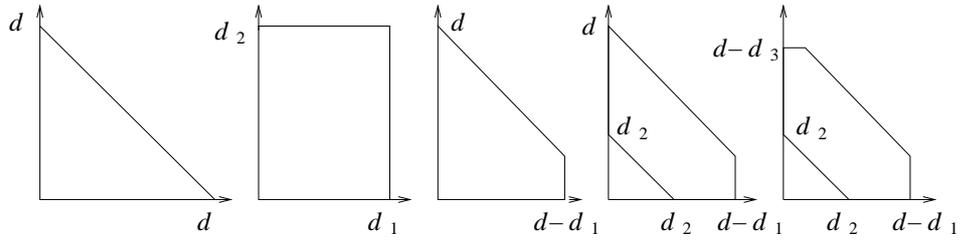}
\end{center}
\caption{Polygons associated with Del Pezzo surfaces}
\label{chf4}
\end{figure}

Let~$\Sig$ be one of these surfaces. An ample divisor on~$\Sig$
defines a linear system, which in suitable toric coordinates is
generated by monomials $x^iy^j$, where $(i,j)$ ranges over all the
integer points of a polygon~$\Del(D)$ of the following form. If
$\Sig = \PP^2$ and $D = d \, \PP^1$, then $\Del(D)$ is the
triangle with vertices $(0,0)$, $(d,0)$, and $(0,d)$. If $\Sig =
\PP^1 \times \PP^1$ and $D$ is of bi-degree $(d_1, d_2)$, then
$\Del(D)$ is the rectangle with vertices $(0, 0)$, $(d_1, 0)$,
$(d_1, d_2)$, and $(0,d_2)$. If $\Sig=\PP^2_k$, $k=1$, $2$, or
$3$, and $D = dL-\sum_{i=1}^k d_iE_i$, then $\Del(D)$ is
respectively the trapeze with vertices $(0,0)$, $(d-d_1,0)$,
$(d-d_1,d_1)$, $(0,d)$, or the pentagon with vertices $(d_2,0)$,
$(d-d_1,0)$, $(d-d_1,d_1)$, $(0,d)$, $(0,d_2)$, or the hexagon
with vertices $(d_2,0)$, $(d-d_1,0)$, $(d-d_1,d_1)$,
$(d_3,d-d_3)$, $(0,d-d_3)$, $(0,d_2)$ (see Figure~\ref{chf4}). The
slopes of the sides of $\Del(D)$ are $0$, $-1$, or $\infty$.

\subsection{Plane tropical curves}\label{plane-curves}
In Sections~\ref{plane-curves} and~\ref{parameterization}
we remind definitions and basic facts
concerning plane tropical curves
(cf.~\cite{Mi}, \cite{GM2}), and fix the notation.

Let $\Del\subset\R^2$ be a nondegenerate convex lattice polygon,
{\it i.e.}, a convex polygon with integer vertices and non-empty
interior. A convex piecewise-linear function
$$F:\R^2\to\R,\quad
F(x)=\max_{\iota\in\Del\cap\Z^2}(\langle\iota,x\rangle+c_\iota),\
\quad\text{where}\quad c_\iota\in\R,\ $$ is called a {\it tropical
polynomial} with Newton polygon $\Del$. Consider the corner locus
$A_F\subset\R^2$ of $F$ ({\it i.e.}, the subset of $\R^2$ where~$F$
is not smooth). The set $A_F$ is naturally stratified and defines a
subdivision~$\Theta_F$ of $\R^2$. The $0$- and $1$-dimensional
elements of the stratification of~$A_F$ are called, respectively,
its {\it vertices} and {\it edges}.

The Legendre transform takes $F$ to a convex piecewise-linear
function $\varrho_F:\Del\to\R$, whose linearity domains $\Del_1$,
$\ldots$, $\Del_N$ are convex lattice subpolygons of $\Del$, and
whose graph is the lower part of the polytope
$\conv\{(\iota,-c_\iota)\in\R^3,\ \iota\in\Del\cap\Z^2\}.$ The
polygons $\Del_1$, $\ldots$, $\Del_N$ give rise to a subdivision
$S_F$ of $\Del$. This subdivision is dual to~$\Theta_F$ in the
following sense: there is a one-to-one correspondence ${\cal D}$
between the elements of $S_F$ and the elements of $\Theta_F$ such
that
\begin{itemize}
\item ${\cal D}$ sends any vertex of $S_F$
to a 2-cell of $\Theta_F$, any edge of $S_F$ to an edge of
$\Theta_F$, and any polygon of $S_F$ to a vertex of $\Theta_F$;
\item for any edge~$e$ of $S_F$, the edge ${\cal D}(e)$
is orthogonal to~$e$;
\item ${\cal D}$ reverses
the incidence relation.
\end{itemize}
Each edge~$e$ of~$A_F$ can be equipped with a weight~$w(e)$ equal
to the {\it lattice length} ({\it i.e.}, the number of integer
points diminished by~$1$) of the dual edge in~$S_F$. The
stratified set $A_F$ whose edges are equipped with the
corresponding weights is called the {\it tropical curve}
associated with the tropical polynomial~$F$. One says that $A_F$
is a plane tropical curve with Newton polygon~$\Del$. A plane
tropical curve determines its Newton polygon~$\Del$ and the dual
subdivision of~$\Del$ uniquely up to translation.

Notice that the unbounded edges of $A_F$ are dual to the edges of
$S_F$ lying on the boundary $\partial\Del$ of~$\Del$. The unbounded
edges of $A_F$ are called {\it ends} of~$A_F$. The unbounded edges
of $A_F$ dual to the edges of $S_F$ which are contained in a side
$\sig$ of $\Del$ are called $\sig$-{\it ends}.

Any edge of a plane tropical curve~$A_F$ has rational slope, and for
any vertex~$v$ of~$A_F$ one has the balancing condition
$$w(e_1)u(v,e_1) + \ldots + w(e_k)u(v,e_k)=0,$$
where $e_1$, $\ldots$, $e_k$ are the
edges adjacent to~$v$, and $u(v,e_i)$ is the primitive integer
vector starting at~$v$ and directed along $e_i$.

The {\it sum} $A_{F^{(1)}} + \ldots + A_{F^{(n)}}$ of plane tropical
curves $A_{F^{(1)}}$, $\ldots$, $A_{F^{(n)}}$ is the plane tropical
curve defined by the tropical polynomial $F^{(1)} + \ldots +
F^{(n)}$. The underlying set of $A_{F^{(1)}} + \ldots + A_{F^{(n)}}$
is the union of underlying sets of $A_{F^{(1)}}$, $\ldots$,
$A_{F^{(n)}}$, and the weight of any edge of~$A_{F^{(1)}} + \ldots +
A_{F^{(n)}}$ is equal to the sum of the weights of the corresponding
edges of summands. A tropical curve in $\R^2$ is {\it reducible} if
it is the sum of two proper tropical subcurves. A non-reducible
tropical curve in $\R^2$ is called {\it irreducible}. A tropical
curve $A_F$ is {\it nodal}, if any polygon of the dual subdivision
$S_F$ is either triangle, or parallelogram. The {\it number of
double points} of a nodal tropical curve~$A_f$ with Newton
polygon~$\Del$ is the sum of the number of parallelograms in~$S_F$
and the number of integer points which belong to the interior
of~$\Del$ and are not vertices of~$S_F$.

The multi-set of vectors $\{w(e)u(e)\}$, where~$e$ runs over the
ends of a plane tropical curve $A_F$ and each primitive integer
vector $u(e)$ is directed along $e$ to infinity, is called the {\it
degree} of $A_F$.

We extend the plane $\R^2$ up to $\widehat\R^2 = \T\times\R$, where
$\T=\R\cup\{-\infty\}$ is equipped with the topology making $\T$
homeomorphic to $[0,+\infty)$ via a logarithmic map, and
correspondingly extend any tropical curve in $\R^2$ by attaching a
vertex on $L_{-\infty} = \{-\infty\}\times\R$ to any horizontal
negatively directed end of the curve.

\subsection{Parameterizations of plane tropical
curves}\label{parameterization} Let $\overline\Gam$ be a finite
connected graph without divalent vertices, and ${\cal V}$ a
collection of certain univalent vertices of $\overline\Gamma$. Put
$\Gam = \overline\Gam \backslash {\cal V}$. Denote by $\Gam^0_0$
the set of non-univalent vertices of $\Gam$, by $\Gam^1$ the set
of edges of~$\Gam$, and by $\Gam^1_\infty$ the set of edges of
$\Gam$ such that the corresponding edges of $\overline\Gam$
terminate at univalent vertices (we call such edges of~$\Gam$ the
{\it ends}).

A {\it parameterized plane tropical curve} is a triple $(\Gamma,
w, h)$, where and $h:\Gam\to\widehat\R^2$ is a continuous proper
map such that
\begin{itemize}\item for
any $E\in\Gam^1$, the restriction of~$h$ to~$E$ is an embedding
into a straight line with a rational slope,
\item if
$E \in \Gam^1_\infty$ and $E$ terminates at a univalent vertex
$V\in\Gam$, then $h(E)$ lies in a horizontal line and $h(V)\in
L_{-\infty}$,
\item for any $V\in\Gam^0_0$,
the union
$$\bigcup_{E\in\Gamma^1, \ V \in \partial E}h(E)$$
is not contained in a line, and one has the balancing condition
$$\sum_{E\in\Gam^1,\ V\in\partial E}
w(E) u(V,E)=0,$$ where $u(V,E)$ is the primitive integer vector
starting at $h(V)$ and directed along $h(E)$.
\end{itemize}
The ends of~$\Gam$ which are adjacent to univalent vertices
of~$\Gam$ are called {\it left}. Due to the properness of $h$, the
non-left ends of~$\Gam$ are mapped onto half-lines which are not
horizontal negatively directed.

The multi-set of vectors $\{w(E)u(E)\ :\ E\in\Gam^1_\infty\}$,
where each primitive integer vector $u(E)$ is directed along
$h(E)$ to infinity, is called the {\it degree} of $(\Gam, w, h)$.

A parameterized plane tropical curve
$(\Gam, w, h)$ is a {\it parameterization}
of a plane tropical curve~$T$ if
\begin{itemize}
\item the image under~$h$ of any vertex of~$\Gam$
is a vertex of~$T$,
\item the closure of the image under~$h$ of any edge of~$\Gam$
is the closure of a union of edges of~$T$,
\item the weight of any edge~$e$ of~$T$ is equal to
$w(E_1) + \ldots + w(E_n)$, where $E_1$, $\ldots$, $E_n$
are the edges of~$\Gam$ whose images under~$h$ contain~$e$.
\end{itemize}
Any parameterized plane tropical curve is a parameterization
of a unique plane tropical curve.
Notice that if $(\Gam, w, h)$ is a parameterization
of a plane tropical curve~$T$, the curve~$T$ might have
vertices that are not images under~$h$ of vertices of~$\Gam$.

The {\it genus} of a parameterized plane tropical curve
$(\Gam,w,h)$ is the first Betti number $b_1(\Gam)$ of~$\Gam$. If a
plane tropical curve~$T$ is irreducible, the minimal genus of its
parameterizations is called the genus of~$T$ and is denoted by
$g(T)$.

The
degree of an irreducible nodal plane tropical curve~$T$
coincides with the degree of any parameterization of~$T$.
If~$T$ is an irreducible nodal plane tropical curve,
then any minimal genus parameterization $(\Gam, w, h)$ of~$T$
is {\it simple}, that is, any vertex of~$\Gam$
has valency either~$3$, or~$1$.
For any two simple parameterizations $(\Gam, w, h)$ and $(\Gam', w', h')$
of a given irreducible nodal plane tropical curve,
there exists
a homeomorphism
$\varphi:\Gam\to\Gam'$ such that $h=h'\circ\varphi$ and
$w(E)=w'(\varphi(E))$ for any $E\in\Gam^1$.

If~$T$ is a nodal plane tropical curve, then each edge of~$T$ is
contained in a unique irreducible subcurve of~$T$. In particular,
a nodal plane tropical curve is uniquely represented as a sum of
its irreducible subcurves. (Notice that this statement is not true
without the nodality assumption.) Furthermore, any irreducible
subcurve of a nodal plane tropical curve is nodal. The genus
$g(T)$ of a nodal plane tropical curve~$T$ is
$g(T^{(1)})+...+g(T^{(n)})-n+1$, where $T^{(1)},...,T^{(n)}$ are
all the irreducible subcurves of~$T$. If~$T$ is a nodal plane
tropical curve with Newton polygon~$\Del$, then the genus of~$T$
is equal to the difference between the number of vertices of the
dual subdivision~$S_T$ which belong to the interior of~$\Del$ and
the number of parallelograms in~$S_T$. In particular, the sum of
the genus of~$T$ and the number of double points of~$T$ is equal
to the number of interior integer points of~$\Del$.

\section{Tropical Welschinger invariants}\label{chsec3}

\subsection{Tropical Welschinger
invariants of toric surfaces}\label{chsec2} Denote by ${\cal C}$
the semigroup of sequences $\alpha = (\alpha_1, \alpha_2,
\ldots)\in\Z^\infty$ with nonnegative terms and finite $l_1$-norm
$\|\alp\| =\sum_i\alp_i$. Each element of~${\cal C}$ contains only
finitely many non-zero terms, so in the description of concrete
sequences we omit zero terms after the last non-zero one. The only
exception concerns the zero element of~${\cal C}$ (the sequence
with all the terms equal to zero). This element is denoted
by~$(0)$. For an element~$\alpha$ in~$\cal C$, put $J\alpha =
\sum_{i=1}^{\infty}(2i-1)\alpha_i$. Define in~${\cal C}$ the
following natural partial order: if each term of a
sequence~$\alpha$ is greater than or equal to the corresponding
term of a sequence~$\beta$, then we say that~$\alpha$ is greater
than or equal to~$\beta$ and write $\alpha \geq \beta$. For two
elements $\alpha = (\alpha_1, \alpha_2, \ldots)$ and $\beta =
(\beta_1, \beta_2, \ldots)$ of~${\cal C}$ such that $\alpha \geq
\beta$, the sequence $\alpha - \beta =(\alpha_1 - \beta_1,\alpha_2
- \beta_2,\dots)$ is an element of~${\cal C}$.

Let $\Del\subset\R^2$ be a nondegenerate convex lattice polygon,
and~$\sig$ the intersection of~$\Del$ with its left vertical
supporting line. Assume that~$\sig$ is a not a point. In this
case, we say that~$\Del$ is {\it left-nondegenerate}. Pick two
elements~$\alpha$ and~$\beta$ in~$\cal C$ such that
$J\alp+J\bet=|\sig|$, where $|\sig|$ is the lattice length of
$\sig$. Fix an integer~$g$, and put
\begin{equation}r = |\partial\Del|-|\sig|+\|\alp\|+\|\bet\|+g-1\
,\label{che201}
\end{equation}
where $|\partial\Del|$ is the
lattice length of the boundary of~$\Del$.
Assume that $r > \|\alp\|$.

Consider the space $\Omega(\Del, \alp,\bet,g) =
(L_{-\infty})^{||\alp||}\times (\R^2)^{r-||\alp||}$ formed by the
(ordered) {\it configurations} $\bx=(\bx^\flat,\bx^\sharp)$ of $r$
points in $\widehat\R^2$ such that $\bx^\flat = (p_1, \ldots,
p_{\|\alp\|})$ is a sequence of $\|\alp\|$ points on
$L_{-\infty}$, and $\bx^\sharp = (p_{\|\alp\|+1}, \ldots, p_r)$ is
a sequence of $|\partial\Del|-|\sig|+\|\bet\|+g-1$ points in
$\R^2$. For any $\bx\in\Omega(\Del,\alp,\bet,g)$, introduce the
set ${\cal T}(\Del,\alp,\bet,g,\bx)$ of nodal plane tropical
curves $T\subset\widehat\R^2$ satisfying the following conditions:
\begin{itemize}
\item $T$ has $\Del$ as Newton polygon and is of genus~$g$;
\item all the $\sig'$-ends of $T$, where $\sig'\ne\sig$, have weight $1$;
\item the number of $\sig$-ends of $T$ is equal to $\|\alp+\bet\|$, and precisely $\alp_i+\bet_i$ of them have weight
$2i-1$, $i\ge 1$;
\item
any irreducible subcurve of~$T$ has a $\sig$-end,
\item $T$ passes through
all the points of~$\bx$,
and
any point $p_k\in\bx^\flat$
is contained in a $\sig$-end of weight $2i_k - 1$,
where the positive integer $i_k$ is determined by
the inequalities
$\sum_{j<i_k}\alp_j<k\le\sum_{j\le i_k}\alp_j$.
\end{itemize}

The first three conditions completely describe the degree of~$T$,
and further on we denote this degree by $\Del^{\alp,\bet}$. For
any $T\in{\cal T}(\Del,\alp,\bet,g,\bx)$ one has
$$
r=\#\End(T)+g-1,
$$
where $\#\End(T)$ is the number of ends of~$T$. For a generic
$\bx\in \Omega(\Del,\alp,\bet,g)$, the set ${\cal
T}(\Del,\alp,\bet,g,\bx)$ is finite. This assertion can be proved
similarly to the corresponding one in~\cite{Mi}. A more precise
statement and a proof is found in
Section~\ref{section-multi-generic}.

Let $T$ be a curve in ${\cal T}(\Del,\alp,\bet,g,\bx)$. If~$T$
does not have edges of even weight, put $W(T) = (-1)^s$, where~$s$
is the total number of integer points lying in the interior of the
triangles of the subdivision~$S_T$ of~$\Del$. Otherwise, put $W(T)
= 0$. The number~$W(T)$ is called the {\it Welschinger
multiplicity} of~$T$. Note that for a reducible tropical curve
$T\in{\cal T}(\Del,\alp,\bet,g,\bx)$, its Welschinger multiplicity
is the product of the Welschinger multiplicities of all the
irreducible subcurves of~$T$.

The Welschinger multiplicity can be also defined for any simply
parameterized plane tropical curve $(\Gam, w, h)$. Namely, let~$V$
be a vertex in~$\Gam^0_0$. If the weights of all the edges
of~$\Gamma$ that are adjacent to~$V$ are odd, denote by $s(V)$ the
number of interior integer points in the triangle built on the
vectors $w(E_1)u(V,E_1)$ and $w(E_2)u(V,E_2)$, where $E_1,E_2$ is
a pair of edges of $\Gam$ adjacent to~$V$, and put $W(V) =
(-1)^{s(V)}$. Otherwise, put $W(V) = 0$. The number
$\prod_{V\in\Gam^0_0}W(V)$ is called the Welschinger multiplicity
of $(\Gam, w, h)$ and is denoted by $W(\Gam, w, h)$. If $(\Gam, w,
h)$ is a simple parameterization of a nodal plane tropical
curve~$T$, then $W(\Gam, w, h) = W(T)$.

Denote
by ${\cal T}^{\irr}(\Del,\alp,\bet,g,\bx)$ the set
of irreducible curves in ${\cal T}(\Del,\alp,\bet,g,\bx)$,
and put
$$W(\Del,\alp,\bet,g,\bx)=
\sum_{T\in{\cal T}(\Del,\alp,\bet,g,\bx)}W(T)\ ,$$ $$
W^{\irr}(\Del,\alp,\bet,g,\bx)
=\sum_{T\in{\cal T}^{\irr}(\Del,\alp,\bet,g,\bx)}W(T)\ .$$

\begin{theorem}\label{cht1}
The numbers $W(\Del,\alp,\bet,g,\bx)$ and
$W^{\irr}(\Del,\alp,\bet,g,\bx)$ do not depend on the choice of a
generic configuration $\bx\in \Omega(\Del,\alp,\bet,g)$.
\end{theorem}

The word `generic' in the statement of Theorem~\ref{cht1} means
that the configurations are taken in an open dense subset of
$\Omega(\Del,\alp,\bet,g)$. This subset is explicitly described in
Section~\ref{section-multi-generic} (see the definition of
multi-tropically generic configurations).

Due to Theorem~\ref{cht1}, one can skip $\bx$ in the notation of
the above numbers. The numbers $W(\Del,\alp,\bet,g)$ (resp.,
$W^{\irr}(\Del,\alp,\bet,g)$) are called {\it multi-component}
(resp., {\it irreducible}) {\it relative tropical Welschinger
invariants}.

The following statement is a corollary
of Mikhalkin's correspondence theorem~\cite{Mi}.

\begin{theorem}\label{chr1}
{\rm (see~\cite{Mi0}, Theorem 6, and \cite{Sh0}, Proposition
6.1).} Let~$\Sig$ be a toric unnodal Del Pezzo surface equipped
with its tautological real structure, $D \subset \Sig$ an ample
divisor, and $\Del$ is a polygon $SL(2, \Z)$-and-translation
equivalent to the polygon $\Del(D)$ defined in
Section~\ref{ample}. Assume that~$\Del$ is left-nondegenerate.
Then,
$$W^{\irr}(\Del,
(0),(|\ \sig|),0)= W(\Sig, D) \ ,$$ where~$\sig$ is the
intersection of~$\Del$ with its left vertical supporting line.
\end{theorem}

The proof of Theorem \ref{cht1} is given in
Section~\ref{invariancy}. It mainly follows the argument of
\cite{GM1,GM2}, where a description of the first order
degenerations and respective bifurcations of simple
parameterizations of tropical curves in count are given. Notice
that Theorem~\ref{cht1} can also be proved via the study of
non-parameterized plane tropical curves in the spirit of
\cite{Sh0}.

\subsection{Recursive formula
for
multi-component
invariants}\label{results}

Denote by~$\theta_k$ the element in ${\cal C}$ whose $k$-th term
is equal to~$1$ and all the other terms are equal to~$0$. For
$\alp,\alp'\in{\cal C}$, $\alp\le\alp'$, put
$$\left(\begin{matrix}
\alp'\\
\alp\end{matrix}\right)
=\prod_{i=1}^\infty\left(\begin{matrix}\alp'_i\\
\alp_i\end{matrix}\right)\ .$$ Extend the definition of the
multi-component relative tropical Welschinger invariants to the
degenerate case $\Del = \sig$ in the following way. If~$\Del$ is a
point, then put
$$W(\Del,(0),(0),g)=\begin{cases}1,\quad & \text{if} \ g=0,\\
0,\quad & \text{otherwise} .\end{cases}$$ If~$\Del$ is a vertical
segment, then put $W(\Del,\alp,\bet,g)=1$ for $\alp + \bet =
(|\sig|)$, $g = 1 - |\sig|$, and put $W(\Del,\alp,\bet,g)=0$, in
all other cases. A number $W(\Del,\alp, \bet, g)$ such that $\Del
= \sig$ will be referred to as an {\it initial value}.

In addition,
put
$W(\Del,\alp,\bet,g) = 0$ whenever $\Del$ is nondegenerate
and $r \leq \|\alp\|$.

Given a convex lattice polygon $\Del$ and a cooriented straight
line $\overrightarrow{s}$ of slope~$0$, $-1$ or $\infty$, take the
supporting straight line $L_{\overrightarrow{s}}$ of $\Del$ such
that $L_{\overrightarrow{s}}$ is parallel to~$\overrightarrow{s}$,
and $\Del$ is contained in the half-plane defined by the
coorientation. Then, the $\overrightarrow{s}$-{\it peeling}
$l_{\overrightarrow{s}}(\Del)$ of $\Del$ is the convex hull of
$(\Del\cap\Z^2)\backslash L_{\overrightarrow{s}}$.

Introduce the set~$\Xi$ formed by the empty set, the lattice
points, the lattice vertical segments, and the convex lattice
left-nondegenerate polygons~$\Del$ such that a primitive integer
normal vector of any face of~$\Del$ belongs to the set $\{(1,0),
(-1,0), (0,1), (0,-1), (1,1), (-1,-1)\}$,

\begin{remark}\label{set-xi}
Any nondegenerate polygon in~$\Xi$ defines a toric unnodal Del
Pezzo surface and an ample divisor on it. The set~$\Xi$ is closed
with respect to the Minkowski sum and any
$\overrightarrow{s}$-peeling $\Del\mapsto
l_{\overrightarrow{s}}(\Del)$, where $\overrightarrow{s}$ is a
cooriented straight vertical line. Furthermore, if the Minkowski
sum of several convex lattice left-nondegenerate polygons is a
polygon in~$\Xi$, then all the summands are in~$\Xi$.
\end{remark}

Let~$\Del \in \Xi$ be
a nondegenerate polygon.
Denote by $\overrightarrow{\infty}$ (respectively,
$\overrightarrow{0}$, $\overrightarrow{-1}$)
a vertical line (respectively, a horizontal line,
a line of slope~$-1$) cooriented by the vector
$(1, 0)$ (respectively, $(0, -1)$, $(1, 1)$).
Denote by $L_{\overrightarrow{0}}$
(respectively, $L_{\overrightarrow{-1}}$)
the support straight line of~$\Del$ such that
$L_{\overrightarrow{0}}$
(respectively, $L_{\overrightarrow{-1}}$)
is parallel to $\overrightarrow{0}$ (respectively, $\overrightarrow{-1}$),
and $\Del$ is contained in the half-plane defined
by the coorientation.
We say that $\Del$ is {\it $\overrightarrow{0}$-nondegenerate}
(respectively, {\it $\overrightarrow{-1}$-nondegenerate}),
if the intersection ${\cal I}$ of~$\Del$ with
$L_{\overrightarrow{0}}$
(respectively, $L_{\overrightarrow{-1}}$) is not a vertex,
and one of the edges neighboring to the edge~${\cal I}$
is of slope $-1$ (respectively, $0$).
Consider a subset~$\daleth$
of $\{\overrightarrow{0}, \overrightarrow{-1}\}$.
The {\it $\daleth$-peeling} $l_\daleth(\Del)$ of~$\Del$ is the result
of the consecutive $\overrightarrow{s}$-peelings
of $l_{\overrightarrow{\infty}}(\Del)$, where $\overrightarrow{s}$ runs
over the elements of~$\daleth$
(note that
in the case $\daleth = \{\overrightarrow{0}, \overrightarrow{-1}\}$,
one has $l_{\overrightarrow{-1}}
(l_{\overrightarrow{0}}(l_{\overrightarrow{\infty}}(\Del)))
= l_{\overrightarrow{0}}
(l_{\overrightarrow{-1}}(l_{\overrightarrow{\infty}}(\Del)))$
since $l_{\overrightarrow{\infty}}(\Del)$ is left-nondegenerate).
The set $\daleth$ is called {\it $\Del$-admissible}
if for any $\overrightarrow{s} \in \daleth$
the polygon $\Del$ is $\overrightarrow{s}$-nondegenerate,
and the polygon $l_\daleth(\Del)$ is either
left-nondegenerate
or a point.
Note
that if $\daleth$ is $\Del$-admissible,
then $l_\daleth(\Del) \in \Xi$.

\begin{theorem}\label{formula}
Let~$\Del \in \Xi$ be a nondegenerate polygon, and~$\sig$ the
intersection of~$\Del$ with its left vertical supporting line.
Then, for any $\alpha,\beta\in{\cal C}$ such that $J\alpha +
J\beta = |\sig|$, and any integer $g$, one has
$$W(\Del,\alpha,\beta,g)=\sum_{\renewcommand{\arraystretch}{0.6}
\begin{array}{c}
\scriptstyle{k\ge 1}\\
\scriptstyle{\beta_k>0}
\end{array}}
W(\Del,\alpha+\theta_k,\beta-\theta_k,g)$$
\begin{equation}+
\sum_{\daleth,\alpha',\beta',g'}
\left(
\begin{matrix}\alpha\\ \alpha'\end{matrix}\right)\left(\begin{matrix}\beta'\\
\beta\end{matrix}
\right)W(l_\daleth(\Del),\alpha',\beta',g')\
,\label{che30}
\end{equation}
where
the latter sum in (\ref{che30}) runs over the
quadruples
$\daleth,\alp',\bet',g'$ satisfying the conditions:

\centerline{$\daleth \subset \{\overrightarrow{0}, \overrightarrow{-1}\}$
is $\Del$-admissible,
$\alp',\bet'\in{\cal C}, \;
g'\in\Z$,}

\vskip-30pt

\begin{equation}
\alp'\le\alp,\quad \bet\le\bet',\quad
(\alp, \bet) \ne (\alp', \bet'),
\quad J\alp' + J\bet'=|\sig'|,\quad g - g' = \|\beta' -
\beta\| - 1\ ,\label{conditions}
\end{equation}
$\sig'$ being the intersection of~$l_\daleth(\Del)$
with its left vertical supporting line.
\end{theorem}

The proof of Theorem~\ref{formula}
basically follows the lines of the proof of Theorem 4.3
from \cite{GM2}, and is presented
in Section~\ref{formula-proof}.

\begin{remark}\label{uni}
The initial values for $W(\Del,\alp,\bet,g)$ {\rm (}{\it i.e.},
the numbers $W(\Del,\alp,\bet,g)$ in the case $\Del = \sig${\rm )}
and the recursive formula given in~{\rm (}\ref{che30}{\rm\/)}
determine all the numbers $W(\Del,\alpha, \beta, g), \Del \in
\Xi$.
\end{remark}

Formula (\ref{che30}) can be seen as a real analogue of the
Caporaso-Harris formula \cite{CH}, Theorem 1.1, and of its
generalizations proposed by R.~Vakil~\cite{Va}. The
Caporaso-Harris formula contains extra coefficients and extra
terms with respect to formula~(\ref{che30}). Comparing two
formulas, one should take into account that a term indexed by
$(\alp, \bet)$ in formula~(\ref{che30}) is an analog of the term
indexed by $(j(\alp),j(\bet))$ of the Caporaso-Harris formula,
where $j: {\cal C} \to {\cal C}$ is the injection associating to a
sequence $\alp = (\alp_1, \alp_2, \alp_3, \ldots) \in {\cal C}$
the sequence $(\alp_1, 0, \alp_2, 0, \alp_3, 0, \ldots)$. In other
words, in formula~(\ref{che30}) we do not consider the tropical
analogs of curves which have even order intersections with the
fixed straight line. Notice also that the Caporaso-Harris formula
contains as a parameter the number of double points instead of the
genus. This difference is not essential, since the genus
determines the number of double points and {\it vice versa}.

\subsection{
Recursive formula for irreducible invariants}\label{chsec4} For
$\alp,\alp^{(1)},...,\alp^{(s)}\in{\cal C}$,
$\alp\ge\alp^{(1)}+...+\alp^{(s)}$, put
$$\left(\begin{matrix}\alp\\ \alp^{(1)},...,\alp^{(s)}\end{matrix}\right)
=\prod_{i=1}^\infty
\frac{\alp_i!}{\alp^{(1)}_i!\cdot...\cdot\alp^{(s)}_i!(\alp_i-\sum_k\alp^{(k)}_i)!}\
.$$ Introduce the set ${\cal S}$ of the $4$-tuples
$(\Del,\alp,\bet,g)$ formed by a polygon $\Del \in \Xi$,
elements~$\alp$ and~$\bet$ in~$\cal C$ such that
$J\alp+J\bet=|\sig|$, where $\sig$ is the intersection of~$\Del$
with its left vertical supporting line, and an integer~$g$. Define
in ${\cal S}$ the following operation:
$$(\Del,\alp,\bet,g)+(\widetilde\Del,\widetilde\alp,\widetilde\bet,\widetilde
g)
=(\Del+\widetilde\Del,\alp+\widetilde\alp,\bet+\widetilde\bet,g+\widetilde
g-1)\ .$$

We extend the definition of the irreducible relative tropical
Welschinger invariants to the degenerate case $\Del = \sig$ in the
following way. If $\Del = \sig$, put $W^\irr(\Del,\alp,\bet,g)=1$
for $g = 0$, $\alp + \bet = (|\sig|)$, and $|\sig| \leq 1$, and
put $W^\irr(\Del,\alp,\bet,g)=0$, in all other cases. A number
$W^{\irr}(\Del,\alp, \bet, g)$ such that $\Del = \sig$ will be
referred to as an {\it initial value}.

In addition,
put
$W^\irr(\Del,\alp,\bet,g) = 0$ whenever $g < 0$.

The irreducible relative tropical Welschinger invariants satisfy a
recursive formula which is similar to the Caporaso-Harris formula
for irreducible relative Gromov-Witten invariants (see~\cite{CH},
Section 1.4).

\begin{theorem}\label{cht3}
Let~$\Del \in \Xi$ and $\sig$
be
as in Theorem~\ref{formula}. Then, for any $\alpha,\beta\in{\cal
C}$ such that $J\alpha + J\beta = |\sig|$, and any integer $g\ge
0$, one has
\begin{equation}
W^\irr(\Del,\alpha,\beta,g)=\sum_{\renewcommand{\arraystretch}{0.6}
\begin{array}{c}
\scriptstyle{k\ge 1}\\
\scriptstyle{\beta_k>0}
\end{array}}
W^\irr(\Del,\alpha+\theta_k,\beta-\theta_k,g)$$
$$+\sum\left(\begin{matrix}\alp\\
\alp^{(1)},...,\alp^{(m)}\end{matrix}\right)
\frac{n!}{n_1!...n_m!} \prod_{i=1}^m\left(
\left(\begin{matrix}\bet^{(i)}\\
\widetilde\bet^{(i)}\end{matrix}\right)
W^\irr(\Del^{(i)},\alp^{(i)},\bet^{(i)},g^{(i)}) \right)\
,\label{che33}\end{equation} where
$$n=|\partial\Del|-|\sig|+\|\bet\|+g-2,\quad
n_i=|\partial\Del^{(i)}|-|\sig^{(i)}|+\|\bet^{(i)}\|+g^{(i)}-1,\
i=1,...,m,$$ and the latter sum in~(\ref{che33})
is taken
\begin{itemize}
\item over all $\Del$-admissible
sets $\daleth \subset \{\overrightarrow{0}, \overrightarrow{-1}\}$,

\item over all
splittings
$$(l_\daleth
(\Del),\alp',\bet',g')=\sum_{i=1}^m(\Del^{(i)},\alp^{(i)},
\bet^{(i)},g^{(i)})$$
in ${\cal S}$, where
$$\displaylines{\alp',\bet'\in{\cal C}, \quad
g'\in\Z, \quad \alp'\le\alp,\quad \bet\le\bet', \cr
J\alp' + J\bet' = |\sig'|, \quad
g - g' = \|\beta' - \beta\| - 1\ ,}$$
$\sig'$ being the intersection of~$l_\daleth(\Del)$
with its left vertical supporting line,
\item over all splittings
$$\bet'=\bet+\sum_{i=1}^m\widetilde\bet^{(i)},\quad\|\widetilde\bet^{(i)}\|>0,\
i=1,...,m\ ,$$ satisfying the restriction
$$\bet^{(i)}\ge\widetilde\bet^{(i)},\ i=1,...,m\ ,$$
\end{itemize}
and factorized by simultaneous
permutations
in the both splittings.
\end{theorem}

\begin{remark}
In the case $g = 0$, the right-hand side of formula~{\rm
(}\ref{che33}{\rm )} reduces to the terms with $g^{(i)} = 0$ and
$\| \widetilde{\bet}^{(i)} \| = 1$.
\end{remark}

The proof of Theorem~\ref{cht3} is a slight modification of the
proof of Theorem~\ref{formula}, and we indicate this modification
at the end of Section~\ref{formula-proof}.

The initial values for $W^\irr(\Del,\alp,\bet,g)$ and
formula~(\ref{che33}) determine all the numbers
$W^\irr(\Del,\alp,\bet,g)$.

\begin{corollary}\label{chc3}
Let~$\Del_1$ and~$\Del_2$ be two nondegenerate polygons in~$\Xi$
such that $\Del_2 \subset \Del_1$. Denote by $\Sig_i$ and $D_i$,
$i=1, 2$, the toric Del Pezzo surface
equipped with its tautological real structure
and the ample divisor on
$\Sig_i$
which are defined by $\Del_i$. Then,
\begin{equation}
W(\Sig_1,D_1)\ge W(\Sig_2,D_2)\ .\label{che70}
\end{equation}
If, in addition, the number of interior integer points of~$\Del_1$
is greater than the number of interior integer points of~$\Del_2$
{\rm (}{\it i.e.}, the genus of a generic member of the linear
system~$|D_1|$ is greater than the genus of a generic member
of~$|D_2|${\rm )}, then
\begin{equation}W(\Sig_1,D_1)>W(\Sig_2,D_2)\
.\label{che71}\end{equation}
\end{corollary}

\begin{corollary}{\rm (cf.  \cite{IKS}, Theorem 1.3)}\label{chc1}
Let~$\Sig$ be a toric unnodal Del Pezzo surface
equipped with its
tautological real structure,
and $D \subset \Sig$ an ample divisor. Then, $W(\Sig, D) > 0$.
\proofend
\end{corollary}

{\bf Proof of Corollaries~\ref{chc3} and~\ref{chc1}}. Since
$\Del_2 \subset \Del_1$, the polygon $\Del_2$ can be obtained from
$\Del_1$ by a sequence of peelings. Thus, it is sufficient to
treat the case when $\Del_2$ is the result of an
$\overrightarrow{s}$-peeling of $\Del_1$. Without loss of
generality, we can assume that $\overrightarrow{s}$ is the left
vertical supporting line of~$\Del_2$ cooriented so that~$\Del_2$
is contained in the half-plane defined by the coorientation.
Denote by~$\sig'$ the intersection of $\overrightarrow{s}$
with~$\Del_2$, and by~$\sig$ the intersection of~$\Del_1$ with its
left vertical supporting line. If~$\sig$ is a point, then
$W(\Sig_1, D_1) = W(\Sig_2, D_2)$. Assume that~$\sig$ is a
nondegenerate segment. One has
$$W(\Sig_1,D_1)=W^\irr(\Del_1,0,(|\sig|),0),\quad
W(\Sig_2,D_2)=W^\irr(\Del_2,0,(|\sig'|),0)\ .$$

The absolute value of the difference
$|\sig|-|\sig'|$ is at most~$1$. According to Theorem~\ref{cht3},
\begin{itemize}
\item if $|\sig'|=|\sig|+1$, then
$$W^\irr(\Del_1,0,(|\sig|),0)\ge|\sig'|\cdot
W^\irr(\Del_2,0,(|\sig'|),0)\ ;$$
\item if $|\sig'|=|\sig|$, then
$$W^\irr(\Del_1,0,(|\sig|),0)\ge
W^\irr(\Del_1,(1),(|\sig|-1),0)$$
$$\ge |\sig'|\cdot W^\irr(\Del_2,0,(|\sig'|),0)
\ ;$$
\item if $|\sig'|=|\sig|-1$, then
$$W^\irr(\Del_1,0,(|\sig|),0)\ge
W^{\irr}(\Del_1,(1),(|\sig|-1),0)$$
$$\ge W^{\irr}(\Del_1,(2),(|\sig|-2),0)
\ge|\sig'|\cdot W^{\irr}(\Del_2,0,(|\sig'|),0)\ .$$
\end{itemize}
Thus, in all the three cases,
\begin{equation}
W^{\irr}(\Del_1,0,(|\sig|),0)\ge|\sig'|\cdot
W^{\irr}(\Del_2,0,(|\sig'|),0)\ .
\label{che74}
\end{equation}
This proves the inequality $W(\Sig_1,D_1)\ge W(\Sig_2,D_2)$.
Since, in addition, $W(\Sig_2, D_2) = 1$ whenever $\Del_2$ does
not have interior integer points, we obtain positivity of the
invariants $W(\Sig, D)$.

If the number of interior integer points of~$\Del_1$ is greater
than the number of interior integer points of~$\Del_2$, then
$|\sig'|\ge 2$, and inequality~(\ref{che74}) implies that
$W(\Sig_1, D_1) > W(\Sig_2, D_2)$. \proofend

\begin{corollary}\label{chc2}
The first six terms of the sequence $W(\PP^2, d\PP^1)$ are as follows:
$$W(\PP^2,\PP^1)=W(\PP^2,2\PP^1)=1, W(\PP^2,3\PP^1)=8,\
W(\PP^2,4\PP^1)=240\ ,$$
$$\hskip3.23cm W(\PP^2,5\PP^1)=18264,\
W(\PP^2,6\PP^1)=2845440\ . \hskip3.23cm \text{\proofend}$$
\end{corollary}

\section{Invariance of tropical Welschinger numbers}\label{invariancy}

\subsection{Moduli spaces of parameterized marked
tropical curves}\label{moduli-spaces} Let $\Del$, $\alp$, $\bet$,
$g$, and $r$ be as in~\ref{chsec2}. A {\it parameterized marked
tropical curve} $(\Gam,w,h,\bpp)$ {\it of degree}
$\Del^{\alp,\bet}$ is a parameterized plane tropical curve
$(\Gam,w,h)$ of degree $\Del^{\alp,\bet}$ equipped with a sequence
$\bpp$ of~$r$ distinct points in~$\Gam$ such that
\begin{itemize}
\item
$\bpp = \bpp^\flat \cup \bpp^\sharp$ starts with
a sequence $\bpp^\flat$ of some univalent
vertices of $\Gam$ and terminates with
a sequence $\bpp^\sharp$ whose points are not
univalent vertices of $\Gam$,
\item the number of points in $\bpp^\flat$ is equal to~$\|\alp\|$,
\item the weight of the ends of $\Gam$ merging to the points
$$P_k\in\bpp^\flat,\qquad\sum_{j<i}\alp_j<k\le
\sum_{j\le i}\alp_j\ ,$$
is $2i-1$, $i\ge 1$,
\item those points
of $\bpp^\sharp$ that coincide with vertices of~$\Gam$ can be
pushed inside the adjacent edges in order to transform $\bpp$ to a
set $\widetilde{\bpp}$ such that the components of
$\overline\Gam\backslash\widetilde{\bpp}$ have no loops and each
of them contains at most one univalent vertex (in particular, the
components of $\overline\Gam\backslash\bpp$ have no loops and each
of them contains at most one univalent vertex).
\end{itemize}
We sometimes call such a sequence~$\bpp$ a {\it configuration}.
The elements of $\bpp$ are called {\it marked points}.

\begin{lemma}\label{one-lemma}
Let~$(\Gam, w, h, \bpp)$ be a parameterized marked tropical curve
such that no point of~$\bpp$ coincides with
a non-univalent vertex of~$\Gam$.
Then, any connected component of $\overline\Gam \setminus \bpp$
contains exactly one univalent vertex.
\end{lemma}

{\bf Proof.}
Identifying all the points of ${\cal V} \subset \overline\Gam$,
we obtain a graph whose first Betti number
is equal to $r - \|\alp\|$. The complement
of~$\bpp^\sharp$ in this graph is a tree,
and the statement follows.
\proofend

Two parameterized marked tropical curves
$(\Gam,w,h,\bpp)$ and
$(\Gam',w',h',\bpp')$
of the same degree
$\Del^{\alp,\bet}$
are called {\it
isomorphic} if there is a homeomorphism
$\varphi:\Gam\to\Gam'$ such that $\varphi(\bpp)=\bpp'$,
and $w(E) = w'(\varphi(E))$ for any $E\in\Gam^1$.
Two
parameterized marked tropical curves
$(\Gam,w,h,\bpp)$ and
$(\Gam',w',h',\bpp')$ of the same degree
$\Del^{\alp,\bet}$
have the same {\it combinatorial type},
if there is a homeomorphism
$\varphi:\Gam\to\Gam'$ such that
\begin{itemize}
\item for any $V \in \Gam^0_0$ and any edge~$E$ adjacent to~$V$,
the vectors $u(V, E)$ and $u(\varphi(V), \varphi(E))$
coincide, where
$u(V, E)$
(respectively, $u(\varphi(V), \varphi(E))$)
is the primitive
integer vector starting at $h(V)$
(respectively, $h'(\varphi(V))$)
and directed along $h(E)$ (respectively, $h'(\varphi(E))$),
\item $w(E) = w'(\varphi(E))$ for any $E\in\Gam^1$,
\item
if a point $P_i \in \bpp$ belongs to
an edge~$E$ (respectively, coincides with a vertex~$V$)
of~$\Gamma$, then the point $P'_i$ belongs
to the edge $\varphi(E)$ (respectively, coincides
with the vertex $\varphi(V)$).
\end{itemize}

Let~$\Lam_{\Del,\alp,\bet,g}$ be the set of all the combinatorial
types of parameterized marked tropical curves $(\Gam, w, h, \bpp)$
of degree $\Del^{\alp,\bet}$ such that $b_1(\Gam) \leq g$ and
$\bpp^\sharp$ contains at least $g - b_1(\Gam)$ points coinciding
with vertices of~$\Gam$. For any $\lam \in
\Lam_{\Del,\alp,\bet,g}$ denote by ${\cal
M}^\lam_{\Del,\alp,\bet,g}$ the set of the isomorphism classes of
parameterized marked tropical curves of combinatorial type~$\lam$.
One can encode the elements $(\Gam, w, h, \bpp)$ of ${\cal
M}^\lam_{\Del,\alp,\bet,g}$ by
\begin{enumerate}
\item[(i)]
the lengths
of images under~$h$
of all the edges
$E\in\Gam^1\backslash\Gam^1_\infty$,
\item[(ii)] the position of
$h(V)\in
\widehat\R^2$ for some vertex $V\in\Gam^0_0$,
\item[(iii)] the coordinates of
the points of $h(\bpp^\flat)$
on $L_{-\infty}$, and
\item[(iv)] the distances of the
points of $h(\bpp^\sharp)$ to the images under~$h$
of certain chosen vertices of the edges
of $\Gam$ which contain the points of $\bpp^\sharp$.
\end{enumerate}
These parameters are called {\it graphic coordinates}. The graphic
coordinates described in the item~(iv) are called {\it marked}.
For a given~$\lam$, graphic coordinates are subject to finitely
many linear equalities and inequalities and identify ${\cal
M}^\lam_{\Del,\alp,\bet,g}$ with the relative interior of a convex
polyhedron in an affine space. We call ${\cal
M}^\lam_{\Del,\alp,\bet,g}$ the {\it moduli space} of
parameterized marked tropical curves of combinatorial type~$\lam$.

A combinatorial type~$\lam' \in \Lam_{\Del,\alp,\bet,g}$ is a {\it
degeneration} of~$\lam \in \Lam_{\Del,\alp,\bet,g}$ if graphic
coordinates on ${\cal M}^\lam_{\Del,\alp,\bet,g}$ and ${\cal
M}^{\lam'}_{\Del,\alp,\bet,g}$ can be chosen in such a way that
${\cal M}^{\lam'}_{\Del,\alp,\bet,g}$ becomes the intersection of
${\cal M}^\lam_{\Del,\alp,\bet,g}$ with some coordinate hyperplanes
in the following sense: the non-zero coordinates of any point of
this intersection are the chosen graphic coordinates of the
corresponding point of ${\cal M}^{\lam'}_{\Del,\alp,\bet,g}$.

For each combinatorial type~$\lambda \in \Lam_{\Del,\alp,\bet,g}$
choose graphic coordinates on ${\cal M}^\lam_{\Del,\alp,\bet,g}$.
Denote by $\overline{{\cal P}^\lam}$ the corresponding polyhedron
and by ${\cal P}^\lam$ its relative interior.
If~$\lam'$ is a degeneration of~$\lam$, then
for any choice of graphic coordinates
on ${\cal M}^\lam_{\Del,\alp,\bet,g}$
there exists a unique choice of graphic coordinates
on ${\cal M}^{\lam'}_{\Del,\alp,\bet,g}$ such that
${\cal M}^{\lam'}_{\Del,\alp,\bet,g}$
becomes the intersection of ${\cal M}^\lam_{\Del,\alp,\bet,g}$
with some coordinate hyperplanes. Denote by $f_{\lam', \lam}$
the affine map identifying this intersection with ${\cal P}^{\lam'}$.

\begin{proposition}\label{boundary}
For any combinatorial type~$\lambda \in \Lam_{\Del,\alp,\bet,g}$
and any face~$F$ of the polyhedron~$\overline{\cal P}^\lam$, there
exists a degeneration~$\lam' \in \Lam_{\Del,\alp,\bet,g}$ such
that $f^{-1}_{\lam',\lam}({\cal P}^{\lam'})=F$.
\end{proposition}

{\bf Proof.} Pick a parameterized marked tropical curve $(\Gam, w,
h, \bpp)$ of combinatorial type~$\lam$. The map~$h$ induces a
metric on~$\Gam$, and thus, an affine structure on it. Any point in ${\cal M}^\lam_{\Del,\alp,\bet,g}$
has a representative $(\Gam, w, \widetilde{h}, \widetilde{\bpp})$
such that $\widetilde{h}$ is affine-linear on any segment that
contains no marked point and no vertex of~$\Gam$. Let~$p \in F$ be
the limit of a sequence of points in ${\cal P}^\lam$. The sequence of
corresponding maps $\widetilde{h}$ converges to a map~$\Gam \to
\widehat\R^2$ which is a composition of a quotient map~$\psi$
from~$\Gam$ to a certain graph~$\Gam'$ and an embedding $h': \Gam'
\to \widehat\R^2$. The sequence of configurations
$\widetilde{\bpp}$ converges to a configuration~$\overline\bpp$ of points
of~$\Gam$. Descending~$w$ and $\overline\bpp$ to~$\Gam'$, one obtains a
parameterized marked tropical curve $(\Gam', w \circ \psi^{-1}, h',
\psi(\overline\bpp))$ together with graphical coordinates identifying it
with~$p$. The combinatorial type
of $(\Gam', w \circ \psi^{-1}, h', \psi(\overline\bpp))$
is a degeneration of~$\lam$.
\proofend

For
any $\lam \in \Lam(\Del, \alp, \bet, g)$,
denote
by ${\cal Q}^\lam$ the projection of ${\cal P}^\lam$
on the coordinate subspace spanned by the non-marked graphical
coordinate axes.

\begin{lemma}{\rm (see~\cite{Mi},
Proposition 2.23)}\label{newdimension}
For any $\lam \in \Lam(\Del, \alp, \bet, g)$,
the dimension of ${\cal Q}^\lam$ is at most~$r$.
Moreover, ${\cal Q}^\lam$ is of dimension~$r$
if and only if the curves in~$\lam$ are
simply parameterized.\proofend
\end{lemma}

If $(\Gam,w,h,\bpp)$ is a simply parameterized marked tropical
curve of degree $\Del^{\alp,\bet}$ such that no point in
$\bpp^\sharp$ is a vertex of~$\Gam$, then the combinatorial type
of $(\Gam,w,h,\bpp)$ is called {\it $\sig$-generic}.

Consider the disjoint union $\coprod_\lam \overline{{\cal P}^\lam}$,
where~$\lam$ runs over all the $\sig$-generic
combinatorial types in $\Lam_{\Del,\alp,\bet,g}$.
Let ${\cal M}_{\Del,\alp,\bet,g}$ be the quotient space
of $\coprod_\lam \overline{{\cal P}^\lam}$ defined by the gluing maps
$f^{-1}_{\lam',\lam_2} \circ f_{\lam',\lam_1}$ for
all the triples $(\lam_1, \lam_2, \lam')$ of combinatorial types
such that $\lam_1$ and $\lam_2$ are $\sig$-generic,
and $\lam'$ is a degeneration of~$\lam_1$ and~$\lam_2$.
We identify
the sets ${\cal M}^\lam_{\Del,\alp,\bet,g}$ with ${\cal P}^\lam$
and consider them
as subspaces in ${\cal M}_{\Del,\alp,\bet,g}$.

In the case $\alp = (0)$, the moduli space ${\cal M}_{\Del,\alp,\bet,g}$
coincides with a moduli space introduced
by Gathmann and Markwig~\cite{GM1}
(in~\cite{GM1} this space is denoted by $\overline{{\cal M}}_{g,\Del}$,
where $\Del = \Del^{(0),\bet}$).

\subsection{Collar}\label{collar}

For a combinatorial type~$\lam \in \Lam_{\Del,\alp,\bet,g}$ of
parameterized marked tropical curves $(\Gam,w,h,\bpp)$ denote
by~$\check\lam$ the combinatorial type of parameterized marked
tropical curves $(\Gam,w,h,\check{\bpp})$ of degree
$\Del^{(0),\alp+\bet}$ such that the configuration $\check{\bpp}$
is obtained from $\bpp$ by pushing the points of $\bpp^\flat$
inside the adjacent edges of~$\Gam$. If~$\lam$ is $\sig$-generic,
so is~$\check\lam$, and {\it vice versa}.

Denote by $\check{\cal M}_{\Del,\alp,\bet,g}$ the union
$\bigcup_{\lam \in \Lam_{\Del,\alp,\bet,g}} {\cal
M}^{\check\lam}_{\Del,(0),\bet,g} \subset {\cal
M}_{\Del,(0),\bet,g}$. Define a map $\Pi: \check{\cal
M}_{\Del,(0),\bet,g} \to {\cal M}_{\Del,\alp,\bet,g}$ associating
to $(\Gam,w,h,\check{\bpp}) \in \check{\cal
M}_{\Del,(0),\bet,g}$ the element $(\Gam,w,h,\bpp)$ in
${\cal M}_{\Del,\alp,\bet,g}$ such that the configuration
$\bpp$ is obtained from $\check{\bpp}$ by moving
each of the first $\|\alp\|$ points of $\check\bpp$
(all these points belong to left ends of~$\Gam$)
to the closest univalent
vertex of~$\Gam$.

Consider the spaces
$$\Omega(\Del, \alp, \bet, g) = (L_{-\infty})^{||\alp||}\times
(\R^2)^{r-||\alp||} \subset (L_{-\infty})^{||\alp||}\times
(\widehat\R^2)^{r-||\alp||}\subset
(\widehat\R^2)^r,$$
the map $\pr: (\widehat\R^2)^r \to
(L_{-\infty})^{||\alp||}\times (\widehat\R^2)^{r-||\alp||}$
replacing the first $||\alp||$ abscissaes by
$-\infty$, and the evaluation map
$$\ev:{\cal M}_{\Del,\alp,\bet,g}\to
(L_{-\infty})^{||\alp||}\times
(\R^2)^{r-||\alp||},\quad
\ev(\Gam,w,h,\bpp)=h(\bpp)\ .$$

\begin{lemma}\label{collar-projection}
For any $\lam \in \Lam_{\Del,\alp,\bet,g}$
the restriction of~$\Pi$ to $\overline{{\cal P}}^{\check\lam}
\cap \check{\cal M}_{\Del,\alp,\bet,g}$
is an affine surjective map to $\overline{{\cal P}}^\lam$,
and its fibers are relative interiors of convex polyhedra
of dimension $||\alp||$. Furthermore,
the restrictions of the evaluation maps to
$\overline{{\cal P}}^\lam$ and $\overline{{\cal P}}^{\check\lam}$
are affine, and $\ev\circ\Pi = \pr\circ\ev$.
\end{lemma}

{\bf Proof}. The restriction of~$\Pi$ to ${\cal P}^{\check\lam}$
is an affine map to ${\cal P}^\lam$. Furthermore, $\Pi$ is
continuous (cf. proof of Proposition~\ref{boundary}). Repeating
the construction of~$\check\lam$ out of $\lam$ one can show that
the restriction of~$\Pi$ to ${\cal P}^{\check\lam}$ is surjective
onto ${\cal P}^\lam$, and each of its fibers is the product
of~$||\alp||$ open rays. The last statement of the lemma is
straightforward. \proofend

\subsection{Zero and one}\label{01}

\begin{lemma}\label{dimension}
For any
combinatorial type $\lam \in \Lam_{\Del,\alp,\bet,g}$,
the dimension of ${\cal M}^\lam_{\Del,\alp,\bet,g}$
is at most $2r - ||\alp||$
{\rm (}recall that $r = |\partial\Del|-|\sig|+\|\alp\|+\|\bet\|+g-1${\rm )},
and $\dim {\cal M}^\lam_{\Del,\alp,\bet,g} = 2r - ||\alp||$ if
and only if~$\lam$ is $\sig$-generic.
\end{lemma}

{\bf Proof.} The statement immediately follows from
Lemma~\ref{newdimension}. \proofend

\begin{lemma}\label{newlemma1}
For any $\sig$-generic combinatorial type~$\lambda$,
\begin{enumerate}
\item[(i)]
the evaluation map~$\ev$ restricts to a bijection between ${\cal
M}^\lam_{\Del,\alp,\bet,g}$ and the interior of a full dimensional
convex polyhedron in $(L_{-\infty})^{||\alp||}\times
(\R^2)^{r-||\alp||}$,
\item[(ii)] any two parameterized marked tropical curves
of combinatorial type~$\lam$ have the same Welschinger multiplicity.
\end{enumerate}
\end{lemma}

{\bf Proof.} Pick a $\sig$-generic combinatorial type~$\lam \in
\Lam_{\Del,\alp,\bet,g}$. Due to Lemma~\ref{collar-projection},
the restriction~$\ev^\lam$ of~$\ev$ to ${\cal P}^\lam$ is affine.
Furthermore, according to~\cite{GM1}, Proposition~4.2,
(cf.~\cite{Mi}) the restriction~$\ev^{\check\lam}$ of~$\ev$ to
${\cal P}^{\check\lam}$ is injective. The equality $\ev \circ \Pi
= \pr \circ \ev$ (see Lemma~\ref{collar-projection}) and the
injectivity of~$\ev^{\check\lam}$ imply the injectivity
of~$\ev^\lam$, since if a point in $\ev({\cal P}^\lam)$ has two
distinct inverse images $(\Gam_1, w_1, h_1, \bpp_1)$ and $(\Gam_2,
w_2, h_2, \bpp_2)$ under $\ev \circ \Pi$, then modifying~$h_1$
and~$h_2$ (alternatively moving some points of $\bpp_1$ and
$\bpp_2$) one gets two distinct points in ${\cal P}^{\check\lam}$
with the same image under~$\ev$.

The second statement of the lemma immediately follows from the
definition of the Welschinger multiplicity. \proofend

The Welschinger multiplicity of the parameterized marked
tropical curves of $\sig$-generic combinatorial type~$\lam$
is denoted by $W(\lam)$.

Our goal is to study the complement in $(L_{-\infty})^{||\alp||}\times
(\R^2)^{r-||\alp||}$ of
the union $\bigcup_\lam \ev({\cal M}^\lam_{\Del,\alp,\bet,g})$,
where $\lam$ runs over all the
combinatorial types in $\Lam_{\Del,\alp,\bet,g}$ such that
the codimension of $\ev({\cal M}^\lam_{\Del,\alp,\bet,g})$
in $(L_{-\infty})^{||\alp||}\times
(\R^2)^{r-||\alp||}$ is at least~$2$.
A combinatorial type $\lam \in \Lam_{\Del,\alp,\bet,g}$
is called {\it injective codimension} $1$ if
${\cal M}^\lam_{\Del,\alp,\bet,g}$ is of codimension~$1$
in ${\cal M}_{\Del,\alp,\bet,g}$, and the restriction
of~$\ev$ to ${\cal M}^\lam_{\Del,\alp,\bet,g}$ is injective.
The set of combinatorial types $\lam \in \Lam_{\Del,\alp,\bet,g}$
such that
the codimension of $\ev({\cal M}^\lam_{\Del,\alp,\bet,g})$
in $(L_{-\infty})^{||\alp||}\times
(\R^2)^{r-||\alp||}$ is at most~$1$ consists of
\begin{itemize}
\item all $\sig$-generic
combinatorial types in $\Lam_{\Del,\alp,\bet,g}$
(see Lemmas~\ref{dimension} and~\ref{newlemma1}) and
\item all the injective codimension~$1$ combinatorial types
in $\Lam_{\Del,\alp,\bet,g}$.
\end{itemize}

\begin{lemma}\label{codim1}
Let $\lam \in \Lam_{\Del,\alp,\bet,g}$
be injective codimension~$1$. Then,
$\check\lam$ is injective codimension~$1$.
Furthermore, any element
in ${\cal M}^\lam_{\Del,\alp,\bet,g}$ is represented by
a parameterized marked tropical curve
$(\Gam,w,h,\bpp)$ such that $b_1(\Gam) = g$ and
\begin{enumerate}
\item[(i)] either $\Gam$ is trivalent, and
exactly one point of $\bpp^\sharp$ is a vertex of $\Gam$,
\item[(ii)] or
one of the vertices of~$\Gam$ is four-valent, the other
non-univalent vertices of~$\Gam$ are trivalent, and
the points of $\bpp^\sharp$ are not vertices of~$\Gam$,
\item[(iii)] or $\Gam$ has two four-valent vertices joined by two
edges~$E_1$ and~$E_2$, the other non-univalent vertices of~$\Gam$
are trivalent, and the points of $\bpp^\sharp$ are not vertices of
$\Gam$.
\end{enumerate}
\end{lemma}

{\bf Proof}. Since all the fibers of~$\Pi$ have the same dimension
(see Lemma~\ref{collar-projection}),
the codimension of ${\cal M}^{\check\lam}_{\Del,\alp,\bet,g}
\subset {\check{\cal M}}_{\Del,\alp,\bet,g}$ is equal to~$1$.
The fibers of~$\pr$ have the same dimension
as the fibers of~$\Pi$. Hence, due to $\ev \circ \Pi = \pr \circ \ev$,
the injectivity of
the restriction of~$\ev$ to ${\cal M}^\lam_{\Del,\alp,\bet,g}$
implies the injectivity of the restriction of~$\ev$
to ${\cal M}^{\check\lam}_{\Del,\alp,\bet,g}$.
The second statement of the lemma follows now from~\cite{GM1},
Proposition 3.9 and Remark 3.6.
\proofend

\subsection{First bifurcation}\label{first-bifurcation}

Let~$\lam \in \Lam_{\Del,\alp,\bet,g}$ be an injective
codimension~$1$ combinatorial type, and $(\Gam, w, h, \bpp)$
a parameterized marked tropical curve of combinatorial type~$\lam$.
Assume that $\Gam$ is trivalent, and
exactly one point of $\bpp^\sharp$ is a vertex of $\Gam$.
Denote this vertex by~$V$.
The last property in the definition of parameterized marked tropical
curves implies that exactly two edges adjacent to~$V$
are {\it allowed} in the following sense: pushing the point
which coincides with~$V$ to any of these edges creates neither a loop
in $\overline\Gam \backslash \bpp$
nor a component of $\overline\Gam \backslash \bpp$
with more than one univalent vertex. Denote the two resulting
$\sig$-generic combinatorial types by~$\lam_+$ and~$\lam_-$.
The combinatorial type~$\lam$ is a degeneration of~$\lam_+$
and~$\lam_-$, and no other $\sig$-generic combinatorial type
has~$\lam$ as degeneration.

\begin{lemma}\label{first}
Let~$\lam$, $\lam_+$, and $\lam_-$ be as above. Then, $\ev({\cal
M}^{\lam_+}_{\Del,\alp,\bet,g})$ and $\ev({\cal
M}^{\lam_-}_{\Del,\alp,\bet,g})$ are on opposite sides of
$\ev({\cal M}^\lam_{\Del,\alp,\bet,g})$.
\end{lemma}

{\bf Proof}. According to Lemma~\ref{codim1}, the combinatorial
type~$\check\lam$ is injective codimension~$1$. As is shown
in~\cite{GM1}, case~(c) in the proof of Theorem 4.8, the
images~$\ev({\cal M}^{\check\lam_+}_{\Del,\alp,\bet,d})$
and~$\ev({\cal M}^{\check\lam_-}_{\Del,\alp,\bet,d})$ of~${\cal
M}^{\check\lam_+}_{\Del,\alp,\bet,g}$ and~${\cal
M}^{\check\lam_-}_{\Del,\alp,\bet,g}$ under the evaluation map are
on opposite sides of $\ev({\cal
M}^{\check\lam}_{\Del,\alp,\bet,d})$. Thus, the statement of the
lemma follows from the relation $\ev \circ \Pi = \pr \circ \ev$
and the fact that~$\pr$ is affine. \proofend

\begin{lemma}\label{first-Welschinger}
Let~$\lam$, $\lam_+$, and~$\lam_-$ be as above. Then, the
Welschinger multiplicities $W(\lam_+)$ and~$W(\lam_-)$ are equal.
\proofend
\end{lemma}

\subsection{Third bifurcation}\label{third-bifurcation}

Let~$\lam \in \Lam_{\Del,\alp,\bet,g}$ be an injective
codimension~$1$ combinatorial type, and $(\Gam, w, h, \bpp)$ a
parameterized marked tropical curve of combinatorial type~$\lam$.
Assume that $\Gam$ has two four-valent vertices~$V$ and~$V'$
joined by two edges~$E_1$ and~$E_2$, the other non-univalent
vertices of~$\Gam$ are trivalent, and the points of $\bpp^\sharp$
are not vertices of $\Gam$. A $\sig$-generic combinatorial
type~$\widetilde\lam$ is a {\it perturbation} of~$\lam$ if
$\widetilde\lam$ is represented by a parameterized marked tropical
curve $(\widetilde\Gam, \widetilde{w}, \widetilde{h},
\widetilde{\bpp})$ such that the graph~$\widetilde\Gam$ is
obtained from~$\Gam$ replacing each four-valent vertex by two
trivalent ones connected by an edge (denote these edges by~$E_V$
and~$E_{V'}$, respectively), and there exists a continuous map
$\varphi: \widetilde\Gam \to \Gam$ satisfying the following
properties:
\begin{itemize}
\item the image of any vertex of~$\widetilde\Gam$ under~$\varphi$
is a vertex of~$\Gam$,
\item if~$E$ is an edge of~$\widetilde\Gam$
different from~$E_V$ and~$E_{V'}$ and a vertex~$W$
is adjacent to~$E$, then the image~$\varphi(E)$
is an edge of~$\Gam$,
the vectors $u(W, E)$ and $u(\varphi(W), \varphi(E))$
coincide,
\item $\varphi(E_V) = V$ and $\varphi(E_{V'}) = V'$,
\item if~$E$ is an edge of~$\widetilde\Gam$
different from~$E_V$ and~$E_{V'}$, then
$\widetilde{w}(E) = w(\varphi(E))$,
\item $\varphi(\widetilde{\bpp}) = \bpp$.
\end{itemize}

\begin{figure}
\begin{center}
\epsfxsize 125mm \epsfbox{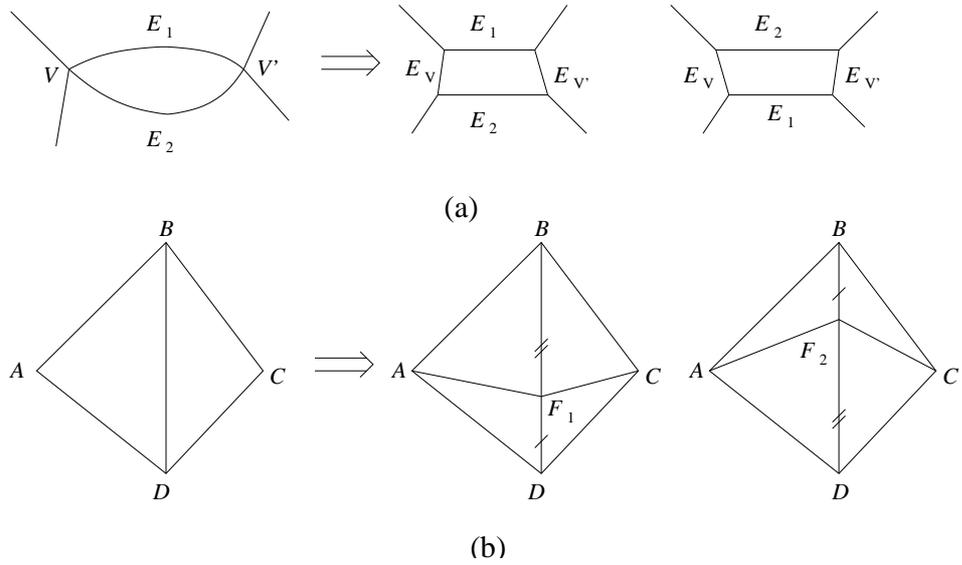}
\end{center}
\caption{Third bifurcation} \label{chfig1}
\end{figure}

\begin{lemma}\label{third}
Let~$\lam \in \Lam_{\Del,\alp,\bet,g}$ be as above. Then, there
are exactly two $\sig$-generic combinatorial types which admit
$\lam$ as a degeneration. These combinatorial types~$\lam_+$
and~$\lam_-$ are perturbations of~$\lam$. The
perturbations~$\lam_+$ and~$\lam_-$ are not equivalent in the
following sense: there is no homeomorphism of their underlying
graphs~$\widetilde\Gam_+$ and~$\widetilde\Gam_-$ which respect the
maps~$\varphi_+$ and~$\varphi_-$ {\rm (}fragments of the
graphs~$\widetilde\Gam_+$ and~$\widetilde\Gam_-$ are shown in
Figure~\ref{chfig1}{\rm (}a{\rm )}{\rm )}. Moreover, $\ev({\cal
M}^{\lam_+}_{\Del,\alp,\bet,g})$ and $\ev({\cal
M}^{\lam_-}_{\Del,\alp,\bet,g})$ are on opposite sides of
$\ev({\cal M}^\lam_{\Del,\alp,\bet,g})$.
\end{lemma}

{\bf Proof}. According to Lemma~\ref{codim1}, the combinatorial
type~$\check\lam$ is injective codimension~$1$. As is shown
in~\cite{GM1}, case~(d) in the proof of Theorem 4.8, there are
exactly two $\sig$-generic combinatorial types~$\check\lam_+$
and~$\check\lam_-$ which admit $\check\lam$ as a degeneration. The
underlying graphs of parameterized marked tropical curves
representing~$\check\lam_+$ and~$\check\lam_-$ are obtained
from~$\Gam$ replacing each four-valent vertex by two trivalent
ones as shown on Figure~\ref{chfig1}(a), and $\ev({\cal
M}^{\check\lam_+}_{\Del,\alp,\bet,d})$ and $\ev({\cal
M}^{\check\lam_-}_{\Del,\alp,\bet,d})$ are on opposite sides of
$\ev({\cal M}^{\check\lam}_{\Del,\alp,\bet,d})$. Thus, the
statements of the lemma follow from the relation $\ev \circ \Pi =
\pr \circ \ev$ and the fact that~$\pr$ is affine. \proofend

\begin{lemma}\label{third-Welschinger}
Let~$\lam$,
$\lam_+$, and~$\lam_-$ be as in Lemma~\ref{third}.
Then,
$W(\lam_+) =
W(\lam_-)$.
\end{lemma}

{\bf Proof}.
Consider parameterized marked
tropical curves $(\widetilde\Gam_+, \widetilde{w}_+, \widetilde{h}_+, \widetilde{\bpp}_+)$
and $(\widetilde\Gam_-, \widetilde{w}_-, \widetilde{h}_-, \widetilde{\bpp}_-)$
of combinatorial types~$\lam_+$ and~$\lam_-$.
The dual subdivisions of plane tropical curves~$T_+$
and~$T_-$
defined by $(\widetilde\Gam_+, \widetilde{w}_+,
\widetilde{h}_+, \widetilde{\bpp}_+)$ and
$(\widetilde\Gam_-, \widetilde{w}_-,
\widetilde{h}_-, \widetilde{\bpp}_-)$
differ by
the fragments shown in Figure \ref{chfig1}(b).

The lattice lengths $|BF_1|$ and $|F_2D|$ are equal to $w(E_1)$,
and the lattice lengths $|F_1D|$ and $|BF_2|$ are equal to
$w(E_2)$. If at least one of the lattice lengths
$|AB|,|BC|,|CD|,|DA|,|BF_1| = |F_2D|,|F_1D| = |BF_2|,$ is even,
then $W(T_+) = W(T_-) = 0$. Assume that all these lengths are odd.
Then $|F_1F_2|$ is even, and therefore $|AF_1| = |AF_2| \mod 2$
and $|F_1C| = |F_2C| \mod 2$. If at least one of the lattice
lengths $|AF_1|$ and $|F_1C|$ is even, then $W(T_+) = W(T_-) = 0$.
If~$|AF_1|$ and~$|F_1C|$ are odd, then the total number~$s_+$ of
interior integer points in the triangles $ABF_1$, $AF_1D$,
$BCF_1$, $F_1CD$ has the same parity as the number~$s_-$ of
interior integer points in the triangles $ABF_2$, $AF_2D$,
$BCF_2$, $F_2CD$. Hence, $W(T_+)= (-1)^{s_+ - s_-} \ W(T_-) =
W(T_-)$. \proofend

\subsection{Second bifurcation}\label{second-bifurcation}

\subsubsection{Preliminaries}\label{prelim}
Let~$\lam \in \Lam_{\Del,\alp,\bet,g}$ be an injective
codimension~$1$ combinatorial type, and let
$(\Gam, w, h, \bpp)$
be
a parameterized marked tropical curve of combinatorial type~$\lam$.
Assume that one of the vertices of~$\Gam$ is four-valent
(denote this vertex by~$V$), the other
non-univalent vertices of~$\Gam$ are trivalent, and
the points of $\bpp^\sharp$ are not vertices of~$\Gam$.
Denote by $E_1$, $E_2$, $E_3$, and $E_4$ the edges of~$\Gam$
which are adjacent to~$V$,
and denote by $L_1$, $L_2$, $L_3$, and $L_4$
the lines containing $h(E_1)$, $h(E_2)$, $h(E_3)$, and $h(E_4)$,
respectively.

A $\sig$-generic combinatorial type~$\widetilde\lam$ is a {\it
perturbation} of~$\lam$ if $\widetilde\lam$ is represented by a
parameterized marked tropical curve $(\widetilde\Gam,
\widetilde{w}, \widetilde{h}, \widetilde{\bpp})$ such that the
graph~$\widetilde\Gam$ is obtained from~$\Gam$ replacing the
vertex~$V$ by two trivalent ones connected by an edge (denote this
edge by~$E_V$), and there exists a continuous map $\varphi:
\widetilde\Gam \to \Gam$ satisfying the following properties:
\begin{itemize}
\item the image of any vertex of~$\widetilde\Gam$ under~$\varphi$
is a vertex of~$\Gam$,
\item if~$E$ is an edge of~$\widetilde\Gam$
different from~$E_V$ and a vertex~$W$ is adjacent to~$E$, then the
image~$\varphi(E)$ is an edge of~$\Gam$, and the vector $u(W, E)$
coincides with $u(\varphi(W), \varphi(E))$,
\item $\varphi(E_V) = V$,
\item if~$E$ is an edge of~$\widetilde\Gam$
different from~$E_V$, then
$\widetilde{w}(E) = w(\varphi(E))$,
\item $\varphi(\widetilde{\bpp}) = \bpp$.
\end{itemize}

\subsubsection{Non-degenerate case}\label{nd}

\begin{figure}
\begin{center}
\epsfxsize 125mm \epsfbox{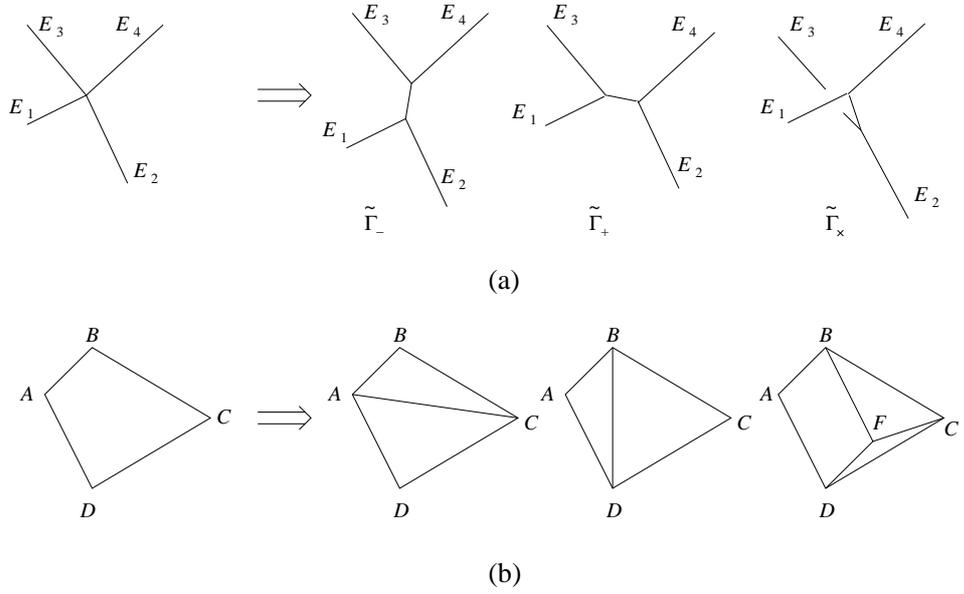}
\end{center}
\caption{Non-degenerate second bifurcation} \label{chfig2}
\end{figure}

\begin{lemma}\label{second1}
Let~$\lam \in \Lam_{\Del,\alp,\bet,g}$ be as in~\ref{prelim}.
Assume that the lines $L_1$, $L_2$, $L_3$, and $L_4$ are pairwise
distinct. Then, there are exactly three $\sig$-generic
combinatorial types which admit $\lam$ as a degeneration. These
combinatorial types~$\lam_\times$, $\lam_+$, and~$\lam_-$ are
perturbations of~$\lam$ and have the following properties:
\begin{itemize}
\item among the edges~$\varphi^{-1}_\times(E_1)$,
$\varphi^{-1}_\times(E_2)$, $\varphi^{-1}_\times(E_3)$,
and $\varphi^{-1}_\times(E_4)$,
there are two edges, $\varphi^{-1}_\times(E_i)$ and
$\varphi^{-1}_\times(E_j)$, such that
the images under $\widetilde{h}_\times$ of their interiors
have a common point,
\item the edges $\varphi^{-1}_+(E_i)$ and $\varphi^{-1}_+(E_j)$
have a common vertex,
\item the edges $\varphi^{-1}_-(E_i)$ and $\varphi^{-1}_-(E_j)$
do not have a common vertex.
\end{itemize}
Fragments of the graphs~$\widetilde\Gam_\times$,
$\widetilde\Gam_+$, and~$\widetilde\Gam_-$
are shown in Figure~\ref{chfig2}{\rm (}a{\rm )}.
The perturbations~$\lam_\times$, $\lam_+$, and~$\lam_-$ are not
equivalent in the same sense as in Lemma~\ref{third}.
Furthermore,
$\ev({\cal M}^{\lam_\times}_{\Del,\alp,\bet,g})$
and $\ev({\cal M}^{\lam_+}_{\Del,\alp,\bet,g})$ are on
the same side of $\ev({\cal M}^\lam_{\Del,\alp,\bet,g})$,
and $\ev({\cal M}^{\lam^-}_{\Del,\alp,\bet,g})$ is on the opposite side.
\end{lemma}

{\bf Proof}. The same arguments as in the proofs of
Lemmas~\ref{first} and~\ref{third} allow one to deduce all the
statements of the lemma from~\cite{GM1}, case~(a) in the proof of
Theorem 4.8. To complete the proof of the last statement, consider
the dual subdivisions of the tropical curves defined by
$(\widetilde\Gam_\times, \widetilde{w}_\times,
\widetilde{h}_\times, \widetilde{\bpp}_\times)$,
$(\widetilde\Gam_+, \widetilde{w}_+, \widetilde{h}_+,
\widetilde{\bpp}_+)$, and $(\widetilde\Gam_-, \widetilde{w}_-,
\widetilde{h}_-, \widetilde{\bpp}_-)$
(see~Figure~\ref{chfig2}(b)), and notice that
\begin{eqnarray}
\Area(BCF)\cdot\Area(CDF) \nonumber
+ \Area(ABD)\cdot\Area(BCD)\\
=
\Area(ABC)\cdot\Area(ACD).\label{che213}\end{eqnarray}
\proofend

\begin{lemma}\label{nondeg-Welschinger}
Let the combinatorial types~$\lam$,
$\lam_\times$, $\lam_+$, and~$\lam_-$ be as in Lemma~\ref{second1}.
Then, $W(\lam_\times) + W(\lam_+) =
W(\lam_-)$.
\end{lemma}

{\bf Proof}. Consider parameterized marked tropical curves
$(\widetilde\Gam_\times, \widetilde{w}_\times,
\widetilde{h}_\times, \widetilde{\bpp}_\times)$,
$(\widetilde\Gam_+, \widetilde{w}_+, \widetilde{h}_+,
\widetilde{\bpp}_+)$, and $(\widetilde\Gam_-, \widetilde{w}_-,
\widetilde{h}_-, \widetilde{\bpp}_-)$ of combinatorial
types~$\lam_\times$, $\lam_+$, and~$\lam_-$. The dual subdivisions
of plane tropical curves~$T_\times$, $T_+$, and~$T_-$ defined by
$(\widetilde\Gam_\times, \widetilde{w}_\times,
\widetilde{h}_\times, \widetilde{\bpp}_\times)$,
$(\widetilde\Gam_+, \widetilde{w}_+, \widetilde{h}_+,
\widetilde{\bpp}_+)$, and $(\widetilde\Gam_-, \widetilde{w}_-,
\widetilde{h}_-, \widetilde{\bpp}_-)$ differ by fragments shown in
Figure~\ref{chfig2}(b). Denote by $W_\times$, $W_+$, and $W_-$ the
(multiplicative) contributions of these fragments to
$W(T_\times)$, $W(T_+)$, and $W(T_-)$. We have to check that
$W_\times + W_+ = W_-$.

If either at least one of the lengths $|AB|,|BC|,|CD|,|AD|$ is
even, or all the lengths $|AC|,|BD|,|FC|$ are even, then
$W_\times=W_+=W_-=0$. So, assume that $|AB|,|BC|,|CD|,|AD|$ are odd,
and $|AC|,|BD|,|FC|$ are not all even.
The first assumption yields
$$\displaylines{
\Area(BCF) = \Area(CDF) = |CF| \mod 2, \cr \Area(ABD) = \Area(BCD)
= |BD| \mod 2, \cr \Area(ABC) = \Area(ACD) = |AC| \mod 2.}$$
Therefore, (\ref{che213}) and the second assumption imply that two
of the lengths $|AC|,|BD|,|FC|$ are odd and the third one is even.
Denote by $N_\times$ (respectively, $N_+$, $N_-$) the total number
of integer points lying in the interior of the triangles $BCF$ and
$CDF$ (respectively, $ABD$ and $BCD$, $ABC$ and $ACD$), and denote
by~$N$ (respectively, $N'$) the number of integer points lying in
the interior of the quadrilateral $ABCD$ (respectively, the
triangle $ABD$). Then, $N_+ = N - |BD| + 1$, $N_- = N - |AC| + 1$,
and $N_\times = N - 2N' - |BD| + 1 - |BF| - |FD| - |CF| + 2$.

If~$|AC|$ is even, while $|BD|$ and $|FC|$ are odd,
then $W_\times=(-1)^{N_\times}$, $W_+=(-1)^{N_+}$, and
$W_-=0$. Furthermore, in this case,
$N_+ = N \mod 2$ and $N_\times = N + 1 \mod 2$,
which yields $W_\times+W_+=0=W_-$.

If~$|BD|$ is even, while $|AC|$ and $|FC|$ are odd,
then $W_\times=(-1)^{N_\times}$, $W_+=0$, and
$W_-=(-1)^{N_-}$. Furthermore, in this case,
$N_- = N \mod 2$ and $N_\times = N \mod 2$,
which yields $W_\times+W_+=(-1)^N=W_-$.

If~$|FC|$ is even, while $|AC|$ and $|BD|$ are odd, then
$W_\times=0$, $W_+=(-1)^{N_+}$, and $W_-=(-1)^{N_-}$. Furthermore,
in this case, $N_+ = N \mod 2$ and $N_- = N \mod 2$, which yields
$W_\times+W_+=(-1)^N=W_-$. \proofend

\subsubsection{Degenerate case}\label{degen}

\begin{figure}
\begin{center}
\epsfxsize 125mm \epsfbox{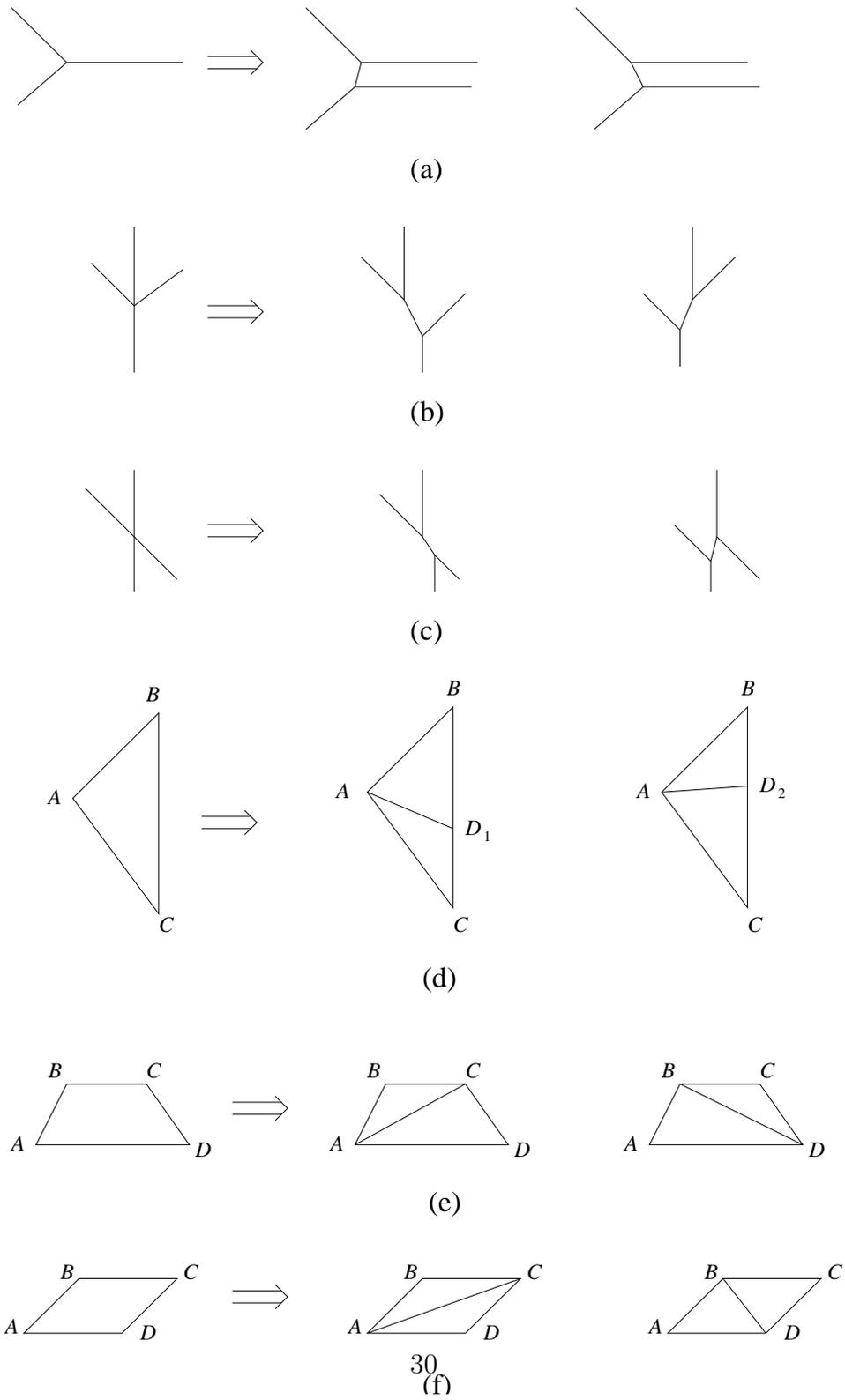}
\end{center}
\caption{Degenerate second bifurcations}\label{chfig3}
\end{figure}

\begin{lemma}\label{second2}
Let~$\lam \in \Lam_{\Del,\alp,\bet,g}$ be as
in~\ref{prelim}.
Assume that some of the four lines
$L_1$, $L_2$, $L_3$, and $L_4$
coincide.
Then, there are exactly two $\sig$-generic combinatorial
types which admit $\lam$ as a degeneration.
These combinatorial types~$\lam_+$ and~$\lam_-$
are perturbations of~$\lam$.
Fragments of the graphs~$\widetilde\Gam_+$ and~$\widetilde\Gam_-$
and their images under~$\widetilde{h}_+$ and $\widetilde{h}_-$
are shown in Figure~\ref{chfig3}{\rm (}a,b,c{\rm )}
{\rm (}in the cases~{\rm (}a{\rm )}, {\rm (}b{\rm )}, and {\rm (}c{\rm )},
the polygon
corresponding to~$V$ in the subdivision dual
to the plane tropical curve defined by $(\Gam, w, h)$
is a triangle,
trapeze, and parallelogram, respectively{\rm )}.
The perturbations~$\lam_+$, and~$\lam_-$ are not
equivalent in the same sense as in Lemma~\ref{third}.
Furthermore,
$\ev({\cal M}^{\lam_+}_{\Del,\alp,\bet,g})$ and
$\ev({\cal M}^{\lam_-}_{\Del,\alp,\bet,g})$
are on the opposite sides
of $\ev({\cal M}^\lam_{\Del,\alp,\bet,g})$.
\end{lemma}

{\bf Proof.} Let $(\Gam, w, h, \bpp)$ be a parameterized marked
tropical curve of combinatorial type~$\lam$. According to
Lemma~\ref{one-lemma} any connected component of $\overline\Gam
\setminus \bpp$ contains exactly one univalent vertex. Assume that
the edges~$E_1$, $E_2$, $E_3$, and $E_4$ are numbered in such a
way that the simple path~$\gam$ in~$\overline\Gam \setminus \bpp$
connecting~$V$ with a univalent vertex contains the edge~$E_4$.

Consider a parameterized marked tropical curve $(\Gam_\circ,
w_\circ, h_\circ, \bpp_\circ)$, where~$\Gam_\circ$ is
obtained from~$\Gam$ by removing the vertex~$V$, the map~$h_\circ$
is obtained by a modification of~$h$ on the edges adjacent to~$V$
in such a way that directions of the images of these edges do not
change, $w_\circ$ is inherited from~$\Gam$, all the points
of~$\bpp_\circ$ but one are inherited from~$\Gam$, and
the additional point~$P_\add$ of~$\bpp_\circ$ belongs
to~$E_4$. The combinatorial type of $(\Gam_\circ, w_\circ,
h_\circ, \bpp_\circ)$ is $\sig$-generic. Denote this
combinatorial type by~$\lam_\circ$. Put $\bx =
h_\circ(\bpp_\circ \backslash \{P_\add\})$ and $p_\add =
h_\circ(P_\add)$. According to Lemma~\ref{newlemma1}, any
two-parameter small perturbation $(\bx(t), p_\add(\tau))$
of $(\bx(0), p_\add(0)) = (\bx, p_\add)$ lifts to a unique continuous
family~$F^{t, \tau} = (\Gam^{t, \tau}_\circ, w^{t, \tau}_\circ,
h^{t, \tau}_\circ, \bpp^{t, \tau}_\circ)$ in the moduli
space of parameterized marked tropical curves of combinatorial
type~$\lam_\circ$. For any  fixed~$t$, the one parameter family
$F^{t, \cdot}$ has the following property: all the edges of $h^{t,
\cdot}_\circ(\Gam^{t, \cdot}_\circ)$ that are not contained
in~$\gam^{t, \cdot}$ preserve their supporting lines.

Among the lines $L_1$, $L_2$, and $L_3$ choose a line $L_i$ such
that the two other lines are distinct, and $L_i$ either coincides
with one of these two lines or coincides with $L_4$. Change, if
necessary, the numbering of lines $L_1$, $L_2$, and $L_3$ in order
to have $i=3$ and the lines $L_2$ and $L_3$ non coinciding. The
position of the lines $L_1^{t, \tau}$, $L_2^{t, \tau}$, and
$L_3^{t, \tau}$ does not depend on~$\tau$, and if $\bx(t) \not\in
\ev({\cal M}^\lam(\Del, \alp, \bet, g))$, then these lines do not
have a common point. Thus, if the perturbation $(\bx(t),
p_\add(\tau))$ is linear, $\bx(t)\not\in \ev({\cal M}^\lam(\Del,
\alp, \bet, g))$ for $t \ne 0$, and $p_\add(\tau) \not\in L_4$ for
$\tau \ne 0$, then $L_3^{t, \tau}$ transports with a non-zero
velocity vector relative to $L_1^{t, \tau} \cup L_2^{t, \tau}$,
while $L_4^{t, \tau}$ (which depends only on $\tau$) performs
another, independent, parallel transport movement with a non-zero
velocity vector relative to $L_1^{t, \tau} \cup L_2^{t, \tau}$.

For a certain sign of~$t$ (assume that this sign is~$-$) and
sufficiently small absolute value of $t$, the ray starting at the
point $L_2^{t, \tau}\cap L_3^{t, \tau}$ and going in the direction
determined by the balancing condition intersects the ray of $F^{t,
\tau} = (\Gam^{t, \tau}_\circ, w^{t, \tau}_\circ, h^{t,
\tau}_\circ, \bpp^{t, \tau}_\circ)$ supported by $L_1^{t, \tau}$.
Selecting $\tau=\tau(t)$ in such a way that $L_4^{t, \tau}$ goes
through the above intersection point, gives rise to a tropical
curve whose combinatorial type is a perturbation of~$\lam$.

If~$L_3$ does not coincide with~$L_1$, then for positive and
sufficiently small values of~$t$, the ray starting at the point
$L_1^{t, \tau}\cap L_3^{t, \tau}$ and going in the direction
determined by the balancing condition intersects the ray of $F^{t,
\tau} = (\Gam^{t, \tau}_\circ, w^{t, \tau}_\circ, h^{t,
\tau}_\circ, \bpp^{t, \tau}_\circ)$ supported by $L_2^{t, \tau}$,
and a construction as above gives rise to a a tropical curve whose
combinatorial type is a perturbation of~$\lam$. If~$L_3$ coincides
with~$L_1$, then for positive and sufficiently small values
of~$t$, the ray starting at the point $L_1^{t, \tau}\cap L_2^{t,
\tau}$ and going in the direction determined by the balancing
condition intersects the ray of $F^{t, \tau} = (\Gam^{t,
\tau}_\circ, w^{t, \tau}_\circ, h^{t, \tau}_\circ, \bpp^{t,
\tau}_\circ)$ supported by $L_3^{t, \tau}$, and again a
construction as above gives rise to a a tropical curve whose
combinatorial type is a perturbation of~$\lam$.

Since some of the lines $L_1$, $L_2$, $L_3$, and $L_4$ coincide,
there are at most two combinatorial types that can be
perturbations of~$\lam$. On the other hand, we constructed two
distinct perturbations of~$\lam$ (denote the first one by~$\lam_-$
and the second one by~$\lam_+$). Moreover, the images $\ev({\cal
M}^{\lam_+}_{\Del,\alp,\bet,g})$ and $\ev({\cal
M}^{\lam_-}_{\Del,\alp,\bet,g})$ are on the opposite sides of
$\ev({\cal M}^\lam_{\Del,\alp,\bet,g})$. \proofend

\begin{lemma}\label{except-Welschinger}
Let~$\lam, \lam_+$, and $\lam_-$
be as in Lemma \ref{second2}.
Then, $W(\lam_+)=W(\lam_-)$.
\end{lemma}

{\bf Proof.} Consider parameterized marked tropical curves
$(\widetilde\Gam_+, \widetilde{w}_+, \widetilde{h}_+,
\widetilde{\bpp}_+)$ and $(\widetilde\Gam_-, \widetilde{w}_-,
\widetilde{h}_-, \widetilde{\bpp}_-)$ of combinatorial
types~$\lam_+$ and~$\lam_-$. The dual subdivisions of plane
tropical curves~$T_+$ and~$T_-$ defined by $(\widetilde\Gam_+,
\widetilde{w}_+, \widetilde{h}_+, \widetilde{\bpp}_+)$ and
$(\widetilde\Gam_-, \widetilde{w}_-, \widetilde{h}_-,
\widetilde{\bpp}_-)$ differ by fragments shown in
Figure~\ref{chfig3}(d), \ref{chfig3}(e), or \ref{chfig3}(f). These
fragments correspond to two splittings of the polygon dual to
$h(V)$, and their (multiplicative) contributions $W_+$ and $W_-$
to $W(T_+)$ and $W(T_-)$, respectively, are equal in each of the
cases~\ref{chfig3}(d,e,f). Indeed, $W_+$ and $W_-$ both vanish if
\begin{itemize}
\item in the case~\ref{chfig3}(d), at least one of the lengths
$|AB|,|AC|,|BD_1|=|CD_2|,|CD_1|=|BD_2|$ in Figure \ref{chfig3}(d)
is even,
\item in the cases~\ref{chfig3}(e,f), at least one of the lengths
$|AB|,|BC|,|CD|,|AD|$ is even.
\end{itemize}
If the
aforementioned lengths are odd then, in the case~\ref{chfig3}(d)
one has $|AD_1| = |AD_2| \mod 2$, and in the cases~\ref{chfig3}(e,f),
one has
$|AC|=|BD|\mod2$.
This yields $W_+=W_-$.
\proofend

\subsection{Tropically generic configurations}\label{generic}
A configuration $\bx\in \Omega(\Del,\alp,\bet,g)$ is called
$(\Del, \alp, \bet, g)$-{\it generic} (resp., {\it almost} $(\Del,
\alp, \bet, g)$-{\it generic}) if the inverse image
$(\ev)^{-1}(\bx)$ of $\bx$ under the evaluation map $\ev: {\cal
M}_{\Del,\alp,\bet,g} \to \Omega(\Del,\alp,\bet,g)$ consists of
parameterized marked tropical curves of $\sig$-generic (resp., of
$\sig$-generic or injective codimension~$1$) combinatorial types.

We say that a parameterized plane tropical curve $(\Gam, w, h)$ of
degree $\Del^{\alp, \bet}$ and genus~$g$ {\it matches} a
configuration $\bx = \bx^\flat \cup \bx^\sharp \in \Omega(\Del,
\alp, \bet, g)$ if $\bx \subset h(\Gam)$, and any point
$p_k\in\bx^\flat$ is contained in the image of a left end of
weight $2i_k - 1$, where the positive integer $i_k$ is determined
by the inequalities $\sum_{j<i_k}\alp_j<k\le\sum_{j\le
i_k}\alp_j$. A $(\Del, \alp, \bet, g)$-generic configuration $\bx
= \bx^\flat \cup \bx^\sharp\in \Omega(\Del,\alp,\bet,g)$ is called
{\it tropically generic} if any parameterized plane tropical curve
matching~$\bx$ and having the degree $\Del^{\alp, \bet}$ and the
genus $g$ defines a plane tropical curve~$T$ which satisfies the
following properties: $T$ is nodal, and no point in~$\bx^\sharp$
coincides with a vertex of~$T$.

\begin{lemma}\label{dimens}
The complement in $\Omega(\Del,\alp,\bet,g)$ of the subset formed
by the tropically generic configurations is a finite closed
polyhedral complex of positive codimension.
\end{lemma}

{\bf Proof.} The subset formed by the non-tropically generic
configurations in $\Omega(\Del,\alp,\bet,g)$ is a projection of
$\cup_{\lam\in\Lam_{\Del, \alp, \bet, g}}{\cal B}^\lam$, where
${\cal B}^\lam \subset {\cal P}^\lam$ is described by a
disjunction of a finite collection of systems of linear equations
and linear inequalities in graphic coordinates. Since the set
$\Lam_{\Del, \alp, \bet, g}$ is finite, it remains to prove that
the tropically generic configurations form an open dense subset in
$\Omega(\Del,\alp,\bet,g)$.

As it immediately follows from Lemma~\ref{dimension} (and
finiteness of $\Lam_{\Del, \alp, \bet, g}$), the $(\Del, \alp,
\bet, g)$-generic configurations form an open dense subset in
$\Omega(\Del,\alp,\bet,g)$. Denote by~$O$ the set of the $(\Del,
\alp, \bet, g)$-generic configurations $\bx \in \Omega(\Del, \alp,
\bet, g)$ such that any parameterized plane tropical curve which
matches~$\bx$ and has the degree $\Del^{\alp, \bet}$ and the genus
$g$ is simply parameterized. Since for any parameterized plane
tropical curve one can always choose a configuration $\bpp \subset
\Gam$ of marked points to obtain a parameterized marked tropical
curve $(\Gam, w, h, \bpp)$, Lemma~\ref{newdimension} implies that
$O$ is an open dense subset in $\Omega(\Del, \alp, \bet, g)$.

Denote by $C_1$ the set of the configurations $\bx = \bx^\flat \cup \bx^\sharp \in O$
admitting a simply parameterized plane tropical
curve $(\Gam, w, h)$
which matches~$\bx$, has the degree $\Del^{\alp,\bet}$ and
the genus~$g$, and satisfies the following property:
the plane tropical curve defined by $(\Gam, w, h)$
has a vertex at one of the points of~$\bx^\sharp$.
Clearly, $C_1$ is a closed nowhere dense subset of~$O$.

Denote by $C_2$ the set of the configurations $\bx \in O$
admitting a simply parameterized plane tropical curve $(\Gam, w,
h)$ which matches~$\bx$, has the degree $\Del^{\alp,\bet}$ and the
genus $g$, and satisfies the following property: the plane
tropical curve defined by $(\Gam, w, h)$ is not nodal. For a
combinatorial type~$\lam \in \Lam(\Del, \alp, \bet, g)$ of simply
parameterized marked tropical curves, consider the image in ${\cal
Q}^\lam$ of the points of ${\cal P}^\lam$ corresponding to
parameterized marked tropical curves defining non-nodal plane
tropical curves. This image is a closed nowhere dense subset of
${\cal Q}^\lam$ as can be shown using the same arguments as in the
proof of~\cite{Mi}, Proposition 2.23. Hence, $C_2$ is a closed
nowhere dense subset of~$O$. Thus, the tropically generic
configurations form an open dense subset in~$\Omega(\Del, \alp,
\bet, g)$. \proofend

\begin{lemma}\label{good-points}
Let $\bx \in \Omega(\Del, \alp, \bet, g)$ be
a tropically generic configuration.
Then, any curve in ${\cal T}^\irr(\Del, \alp, \bet, g, \bx)$
can be parameterized by an element of the inverse image $(\ev)^{-1}(\bx)$
of~$\bx$ under the evaluation map
${\cal M}_{\Del, \alp, \bet, g} \to \Omega(\Del, \alp, \bet, g)$.
\end{lemma}

{\bf Proof.}
Pick an element~$T$ in ${\cal T}^\irr(\Del, \alp, \bet, g, \bx)$,
and consider a simple parameterization $(\Gam, w, h)$ of~$T$.
The configuration~$\bx$ lifts to
a configuration~$\bpp \subset \Gam$.
Assume that
$\overline \Gam \backslash \bpp$ has a component containing either a loop,
or two univalent vertices.
Then, there exists a one-dimensional family
of simply parameterized plane tropical curves $(\Gam, w, h_t)$
such that $h_t(\bpp) = \bx$ for any~$t$, and the coordinates
of the images of vertices of~$\Gam$ depend linearly on~$t$.
Hence, this family degenerates either to a situation
of collision of two vertices of~$\Gam$, or to a situation
of collision of a point in~$\bpp$ and a vertex of~$\Gam$.
The both cases contradict the fact that $\bx$ is tropically generic.
\proofend

\begin{lemma}\label{finiteness}
For any tropically generic configuration $\bx\in \Omega(\Del,\alp,\bet,g)$,
the set ${\cal T}^\irr(\Del,\alp,\bet,g,\bx)$ is finite.
\end{lemma}

{\bf Proof.}
The lemma follows from Lemma~\ref{newlemma1}, Lemma~\ref{good-points},
and finiteness of the set $\Lam_{\Del, \alp, \bet, g}$.
\proofend

\subsection{Multi-tropically generic configurations}\label{section-multi-generic}
Let ${\cal S}_\st$ be the set of the $4$-tuples $(\Del_\st,
\alp_\st, \bet_\st, g_\st)$ formed by a left-nondegenerate convex
lattice polygon $\Del_\st$, elements $\alp_\st$ and $\bet_\st$ in
$\cal C$ such that
$J\alp_* + J\bet_*=|\sig_\st|$
(where
$\sig_\st$ is the intersection of~$\Del_\st$ with its left vertical
supporting line), and an integer~$g_*$. Define an addition operation
in ${\cal S}_*$ in the same way as in~${\cal S}$ (see
Section~\ref{chsec4}).

A partition
$\bigsqcup_{j=1}^\newl \bx^{(j)}$
of~$\bx \in \Omega(\Del,
\alp, \bet, g)$ is called {\it compatible} with a splitting
$$(\Del,\alp,\bet,g)=\sum_{\newj=1}^\newl(\Del^{(\newj)},
\alp^{(\newj)},\bet^{(\newj)},g^{(\newj)})$$ in ${\cal S}_\st$ if
$\bx^{(\newj)} \in \Omega(\Del^{(\newj)},
\alp^{(\newj)},\bet^{(\newj)},g^{(\newj)})$ for any $\newj = 1$,
$\ldots$, $\newl$. A configuration $\bx\in
\Omega(\Del,\alp,\bet,g)$ is called {\it multi-tropically generic}
if for any splitting
$$(\Del,\alp,\bet,g)=\sum_{\newj=1}^\newl(\Del^{(\newj)},
\alp^{(\newj)},\bet^{(\newj)},g^{(\newj)})$$
in ${\cal S}_\st$ and
any partition $\bx = \bigsqcup_{\newj=1}^\newl \bx^{(\newj)}$
compatible with this splitting
the following holds:
\begin{itemize}
\item
each configuration $\bx^{(\newj)}$ is tropically generic,
\item
any sum $\sum_{\newj = 1}^\newl T^{(\newj)}$ is a nodal plane
tropical curve whenever $T^{(\newj)}$ is a plane tropical curve
defined by a parameterized tropical curve which belongs to the
inverse image $(\ev)^{-1}(\bx^{(\newj)})$ of $\bx^{(\newj)}$ under
the evaluation map $\ev: {\cal M}_{\Del^{(\newj)},
\alp^{(\newj)},\bet^{(\newj)},g^{(\newj)}} \to
\Omega(\Del^{(\newj)},\alp^{(\newj)},\bet^{(\newj)},g^{(\newj)})$.
\end{itemize}

\begin{lemma}\label{multi-generic}
The complement in $\Omega(\Del,\alp,\bet,g)$ of the subset formed
by the multi-tropically
generic configurations is a finite closed polyhedral complex
of positive codimension.
\end{lemma}

{\bf Proof.} Consider the subset $U \subset
\Omega(\Del,\alp,\bet,g)$ formed by the configurations $\bx$ such
that, for any splitting
$(\Del,\alp,\bet,g)=\sum_{\newj=1}
^\newl
(\Del^{(\newj)},
\alp^{(\newj)},\bet^{(\newj)}, g^{(\newj)})$ in ${\cal S}_\st$ and
any partition $\bx = \bigsqcup_{\newj=1}^\newl \bx^{(\newj)}$
compatible with the given splitting, all the configurations
$\bx^{(1)}$, $\ldots$, $\bx^{(\newl)}$ are tropically generic.
Lemma~\ref{dimens} implies that $U$ is an open dense subset of
$\Omega(\Del,\alp,\bet,g)$. Now the statement of the lemma follows
from the fact that the sum of several nodal plane tropical curves
can be always made nodal by arbitrarily small parallel shifts of
summands. \proofend

\begin{lemma}\label{multi-finiteness}
For any multi-tropically generic configuration $\bx\in \Omega(\Del,\alp,\bet,g)$,
the set ${\cal T}(\Del,\alp,\bet,g,\bx)$ is finite.
\end{lemma}

{\bf Proof.}
The statement is an immediate consequence of Lemma~\ref{finiteness}.
\proofend

\begin{lemma}\label{chl5}
Let~$\bx \in \Omega(\Del, \alp, \bet, g)$ be a multi-tropically
generic configuration. Then,
\begin{equation}
W(\Del,\alpha,\beta,g,\bx) = \sum
\prod_{\newj=1}^\newl
W^{\irr}(\Del^{(\newj)},
\alpha^{(\newj)},\beta^{(\newj)},g^{(\newj)}, \bx^{(\newj)})\ ,
\label{che34}
\end{equation}
where the sum is taken over all unordered splittings
$(\Del,\alp,\bet,g)=\sum_{\newj=1}^\newl(\Del^{(\newj)},
\alp^{(\newj)},\bet^{(\newj)},g^{(\newj)})$ in ${\cal S}_\st$ and
all compatible partitions of~$\bx$.
\end{lemma}

{\bf Proof.} Straightforward.
\proofend

\subsection{Proof of Theorem~\ref{cht1}}\label{proof2}

Notice that it is
enough to establish the invariance of
the numbers $W^\irr(\Del,\alpha,\beta,g)$, since the invariance of
$W(\Del,\alpha,\beta,g)$ will then follow from~(\ref{che34}).

Pick two tropically generic configurations $\bx$ and $\by$ in
$(L_{-\infty})^{||\alp||}\times (\R^2)^{r-||\alp||}$, and connect
them by a path~$\xi \subset (L_{-\infty})^{||\alp||}\times
(\R^2)^{r-||\alp||}$ such that
\begin{itemize}
\item $\xi$ consists of $(\Del, \alp, \bet, g)$-generic configurations
and finitely many almost $(\Del, \alp, \bet, g)$-generic configurations,
\item for any almost $(\Del, \alp, \bet, g)$-generic
configuration~$\bz \in \xi$
and any injective codimension~$1$ combinatorial type
of parameterized marked tropical curves
in $(\ev)^{-1}(\bz) \subset {\cal M}_{\Del,\alp,\bet,g}$,
the path~$\xi$ intersects $\ev({\cal M}^\lam_{\Del,\alp,\bet,g})$
transversally.
\end{itemize}
According to Lemmas~\ref{first-Welschinger},
\ref{third-Welschinger}, \ref{nondeg-Welschinger},
and~\ref{except-Welschinger}, the value of $\sum_{\lam \ \in \
{\ev}^{-1}(\bz)}W(\lam)$ is the same for all $(\Del, \alp, \bet,
g)$-generic configurations $\bz \in \xi$. In particular,
$W^\irr(\Del,\alp,\bet,g,\bx) = W^\irr(\Del,\alp,\bet,g,\by)$.
\proofend

\section{Proof of
the recursive formulas}\label{formula-proof}

\subsection{Auxiliary lemmas}\label{auxiliary}

\begin{lemma}\label{splitting}\text{\rm(cf.~\cite{GM2}, proof of Theorem 4.3)}.
Let $\Del$, $\alp$, $\bet$, and $g$ be as in
Theorem~\ref{formula}. Fix positive real numbers~$\varepsilon$
and~$N$, and consider a multi-tropically generic configuration
$\bx = (\bx^\flat, \bx^\sharp) \in \Omega(\Del, \alp, \bet, g)$
such that
\begin{itemize}
\item the second coordinates of all the points in~$\bx$
belong to the interval $(-\varepsilon, \varepsilon)$,
\item the first coordinate of one point in~$\bx^\sharp$
is smaller than $-N$,
while the first coordinates of all other points in~$\bx^\sharp$
belong to the interval $(-\varepsilon, \varepsilon)$.
\end{itemize}
Then, for any tropical curve $T \in {\cal T}(\Del,\alpha, \beta,
g, \bx)$, the second coordinates of all trivalent vertices of~$T$ belong
to the interval $(-\varepsilon, \varepsilon)$. Furthermore, if~$N$
is sufficiently
large
with respect to~$\varepsilon$, there
exist real numbers~$a$ and~$b$ satisfying the inequalities
$-N < a < b < -\varepsilon$ and satisfying the
following condition: for any
tropical curve $T \in {\cal T}(\Del,\alpha, \beta, g, \bx)$
the intersection of~$T$ with the rectangle
$\{(x, y) \in \R^2: \; a \leq x \leq b \; \text{\rm and} \;
-\varepsilon \leq y \leq \varepsilon\}$ does not contain vertices
of~$T$ and consists of horizontal segments.
\end{lemma}

{\bf Proof.} Consider a tropical curve~$T \in {\cal T}(\Del,
\alpha, \beta, g, \bx)$. Among the trivalent vertices of~$T$
choose a vertex $v = (v_1, v_2)$ having the maximal second
coordinate. The curve~$T$ has an end starting at~$v$ and pointing
upwards. This end is orthogonal to one of the upper sides
of~$\Del$, and thus, is of weight~$1$ and of direction either $(0,
1)$ or $(1, 1)$. Hence, $T$ should have another edge which starts
at~$v$ and does not point downwards. Consider a simple
parameterization $(\Gam, w, h)$ of the irreducible subcurve $T^0$
of~$T$ such that $v \in T^0$, and denote by~$\bpp$ the lifting
of~$\bx$ to $\Gam$. If $v_2 > \varepsilon$, the connected
component of $\overline\Gam \setminus \bpp$ containing $h^{-1}(v)$
has at least two ends, and thus, $(\Gam, w, h, \bpp)$ is not a
parameterized marked tropical curve. This contradicts
Lemma~\ref{good-points}. In the same way one can show that the
curve~$T$ does not have vertices below the line $y =
-\varepsilon$. This proves the first statement of the lemma.

Denote by~$R$ the rectangle $\{(x, y) \in \R^2: \; -N \leq x \leq
-\varepsilon \; \text{\rm and} \; -\varepsilon \leq y \leq
\varepsilon\}$. Let~$T^1$ be an irreducible subcurve of~$T$ such
that the intersection of~$T^1$ with the interior of~$R$ is
non-empty, and let~$(\Gam^1, w^1, h^1)$ be a simple
parameterization of~$T^1$. As follows from Lemma~\ref{good-points}
and the first statement of the current lemma, the image under~$h$
of any path~$\gamma \subset \overline\Gamma \setminus \bpp$ does
not intersect one of the two horizontal edges of~$R$. Thus, for
any point $(x_1, y_1)$ belonging to the interior of~$R$  and to a
non-horizontal edge of~$T^1$, there exists a path $\gamma \subset
\overline\Gamma \setminus \bpp$ such that $h(\gamma)$
contains~$(x_1, y_1)$ and is the graph of a strictly monotone
function~$f$ defined either on the interval $[-N, x_1]$, or on the
interval $[x_1 -\varepsilon]$. Since there are only finitely many
slopes that can be realized by the edges of a tropical curve with
Newton polygon~$\Del$, the length of the definition interval
of~$f$ is bounded from above by a constant depending only
on~$\Del$. This proves the second statement of the lemma.
\proofend

\begin{lemma}\label{lch26}
Consider a non-degenerate lattice polygon~$\del$ whose projection
to the horizontal coordinate axis coincides with the segment
$[0,1]$. Put $\sig_1=\del\cap\{x=0\}$, $\sig_2=\del\cap\{x=1\}$,
and introduce the vectors $u_1=(-1,0)$, $u_2=(1,0)$. Fix a point
$p \in \R^2 \subset \widehat\R^2$ and an ordered splitting
$|\sig_1|=n_{1,1}+...+n_{1,m_1}$ {\rm (}respectively,
$|\sig_2|=n_{2,1}+...+n_{2,m_1}${\rm )} of $|\sig_1|$ {\rm
(}respectively, of $|\sig_2|${\rm )} into positive integer
summands. Fix also two non-increasing sequences of real numbers
$y_{1,1},...,y_{1,m_1}$ and $y_{2,1},...,y_{2,m_2}$. {\rm (}In the
case $|\sig_i| = 0$, the sequences $n_{i,1}, \ldots , n_{i,m_i}$
and $y_{i,1},...,y_{i,m_i}$ are empty.{\rm )} Then, there exists a
parameterized plane tropical curve $(\Gam,w,h)$ satisfying the
following conditions:
\begin{enumerate}
\item[(1)] $\Gam$ has genus zero and degree
$(n_{1,1}u_1,\ ...\ ,\ n_{1,m_1}u_1,\ n_{2,1}u_2,\ ...\ ,\ n_{2,m_2}u_2,\
u_3,\ u_4)$, where~$u_3$ and~$u_4$ are primitive integer
outward normal vectors of the two non-vertical
sides of~$\del$,
\item[(2)]
$\Gam$ has $m_1$ ends $E_{1,1}$, $\ldots$, $E_{1,m_1}$ such that,
for any $k = 1$, $\ldots$, $m_1$, the end~$E_{1,k}$ is of
weight~$n_{1,k}$, and $h(E_{1,k})$ is a horizontal negatively
directed ray which is contained in the line $y = y_{1,k}$,
\item[(3)]
$\Gam$ has $m_2$ ends $E_{2,1}$, $\ldots$, $E_{2,m_2}$
such that, for any $k = 1$, $\ldots$, $m_2$,
the end~$E_{2,k}$ is of weight~$n_{2,k}$, and
$h(E_{2,k})$ is a horizontal positively directed ray
which is contained in the line $y = y_{2,k}$,
\item[(4)] $p\in h(\Gam)$.
\end{enumerate}
Furthermore, if all the numbers
$y_{1,1},...,y_{1,m_1},y_{2,1},...,y_{2,m_2}$ differ from the
second coordinate of $p$, then all parameterized plane tropical
curves having the above properties define the same plane tropical
curve, and $p$ belongs to the interior of a non-horizontal edge of
this curve.
\end{lemma}

\begin{figure}
\begin{center}
\epsfxsize 125mm \epsfbox{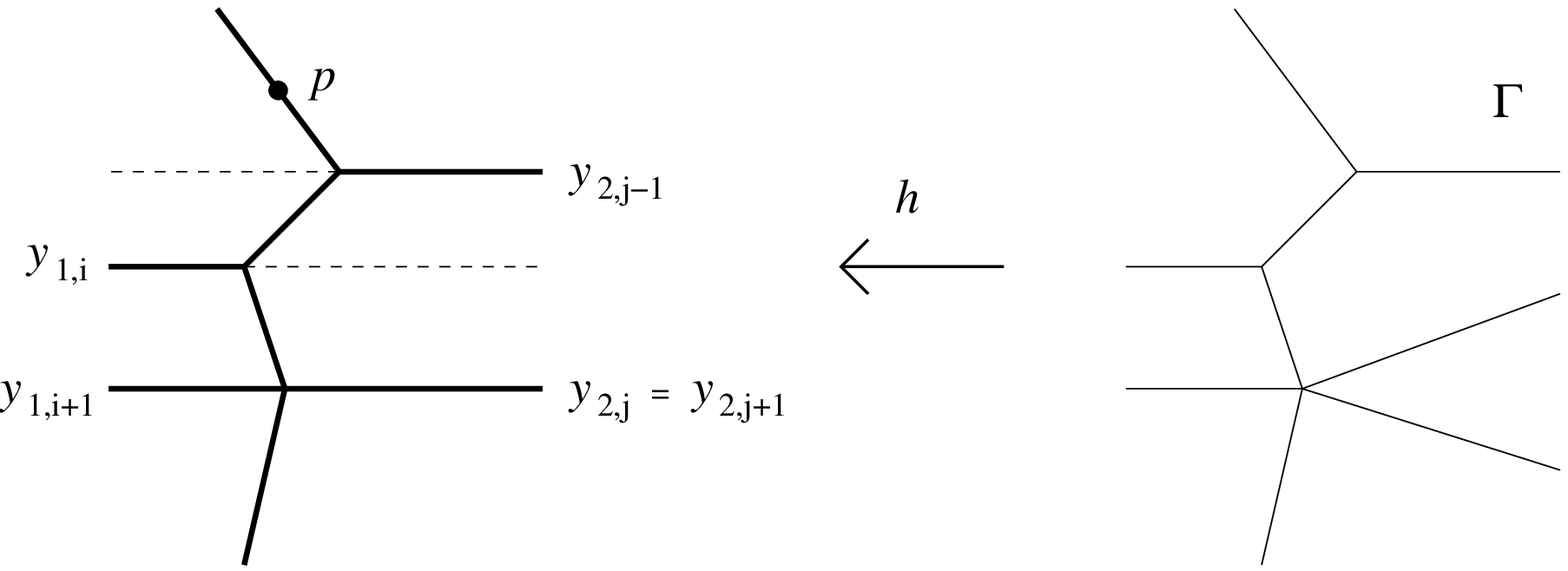}
\end{center}
\caption{ Curve $(\Gam, w, h)$ in Lemma \ref{lch26}}
\label{chfig9}
\end{figure}

{\bf Proof.} Any plane tropical curve having a parameterization
with the described properties can be constructed in the following
way. Take the union $\beth$ of all the lines $y = y_{1,1}, \ldots,
y = y_{1,m_1}, y = y_{2,1}, \ldots, y = y_{2,m_2}$, and consider a
broken line~${\cal L}$ such that
\begin{itemize}
\item any vertex of~$\cal L$ belongs to~$\beth$,
\item any edge of~$\cal L$ has a rational nonzero slope,
\item the two unbounded edges of~$\cal L$ have
directions determined by the vectors~$u_3$ and~$u_4$,
\item for any vertex~$v$ of~$\cal L$,
the two primitive integer vectors~$e_1$ and~$e_2$
starting at~$v$ and directed along the adjacent edges of~$\cal L$
satisfy the relation
$$e_1 + e_2 + \sum_{y_{1,k} = v_y}n_{1,k}u_1 +
\sum_{y_{2,k} = v_y}n_{2,k}u_2 = 0,$$
where $y = v_y$ is the line containing~$v$.
\end{itemize}
Make a horizontal shift of~$\cal L$ in order to obtain a broken
line containing~$p$, and extend the result (in a unique possible
way) to a plane tropical curve which admits a parameterization
satisfying the properties (1) - (3). This proves the existence.
The second statement of the lemma immediately follows from the
construction. \proofend

Consider the subset~${\cal X} \subset \Omega(\Del, \alp, \bet, g)$
formed by the multi-tropically generic
configurations $\bx\ = \bx^\flat \cup \bx^\sharp$
satisfying the following property:
for any point $p\in\bx^\sharp$,
there exists a number~$M$ such that
any configuration
$\widetilde\bx \in \Omega(\Del, \alp, \bet, g)$
obtained from
$\bx$ by a horizontal shift of~$p$
to a point whose first
coordinate is smaller than~$M$ and different from $-\infty$
is multi-tropically generic.

\begin{lemma}\label{lch28}
The subset~${\cal X}$ is open dense in $\Omega(\Del, \alp, \bet,
g)$.
\end{lemma}

{\bf Proof.} Straightforward from Lemma~\ref{multi-generic}.
\proofend

\begin{lemma}\label{lch27}
Let $\Del$,
$\alp$, $\bet$, and $g$ be as in
Theorem~\ref{formula}. Choose $\daleth$, $\alp'$, $\bet'$, and $g'$
satisfying the conditions~{\rm (}\ref{conditions}{\rm )}, a
configuration $\bx = \bx^\flat \cup \bx^\sharp \in {\cal X}$, and
a point~$p \in \bx^\sharp$. Then any subconfiguration $\by \subset
\bx$ such that $\by \in \Omega(l_\daleth(\Del), \alp', \bet', g')$
and $\by^\sharp=\bx^\sharp\backslash\{p\}$ is
multi-tropically generic.
\end{lemma}

{\bf Proof.}
Consider a splitting
$(
l_\daleth(\Del),\alp',\bet',g')=
\sum_{\newi=1}^\newm(\Del^{(\newi)},
\alp^{(\newi)},\bet^{(\newi)}, g^{(\newi)})$ in ${\cal S}_\st$,
and a partition $\bigsqcup_{\newi=1}^\newm \by^{(\newi)}$ of~$\by$
compatible with this splitting. All the polygons $\Del^{(1)}$,
$\ldots$, $\Del^{(\newm)}$ are in~$\Xi$
(it follows from Remark~\ref{set-xi} and the fact
that $l_\daleth(\Del) \in \Xi$).
For each $\newi = 1$, $\ldots$, $m$, pick a parameterized plane
tropical curve $(\Gam^{(\newi)},w^{(\newi)},h^{(\newi)})$
matching~$\by^{(\newi)}$ and having the degree
$(\Del^{(\newi)})^{\alp^{(\newi)},\bet^{(\newi)}}$ and the genus
$g^{(\newi)}$.

We say that a left end of~$\Gam^{(\newi)}$ is {\it marked} if this
end terminates at a marked univalent vertex of~$\Gam^{(\newi)}$.
Among the non-marked left ends of $\coprod_{\newi =
1}^\newm\Gam^{(\newi)}$ choose $||\bet'-\bet||$ ends whose weights
fit the sequence $\bet'-\bet$, and denote the chosen set by~${\cal
E}$. Consider a polygon~$\del$ such that
\begin{itemize}
\item the projection of~$\del$ to the horizontal axis
coincides with the segment $[0, 1]$,
\item the sides $\sig_1=\del\cap\{x=0\}$
and $\sig_2=\del\cap\{x=1\}$ satisfy $|\sig_1|=J(\alp)-J(\alp')$,
$|\sig_2|=J(\bet')-J(\bet)$,
\item the two non-vertical sides
of~$\del$ are respectively parallel to the sides
$\top(\Del, \daleth)$ and $\bot(\Del, \daleth)$ of $\Del$
which are defined
as follows:
the side $\top(\Del, \daleth)$ is of slope~$-1$ and not
a neighbor of~$\sig$
if $\overrightarrow{0} \in \daleth$, and
$\top(\Del, \daleth)$ is the non-vertical side
adjacent to the upper vertex of~$\sig$
otherwise;
the side $\bot(\Del, \daleth)$ is of slope~$0$ and not a neighbor of~$\sig$
if $\overrightarrow{-1} \in \daleth$, and
$\bot(\Del, \daleth)$ is the non-vertical side
adjacent to the lower vertex of~$\sig$
otherwise.
\end{itemize}
The conditions~(\ref{conditions}) imply that~$\del$ is
nondegenerate. Put $m_1 = \|\alp - \alp'\|$, consider a
non-increasing sequence $y_{1,1}$, $\ldots$, $y_{1,m_1}$ of second
coordinates of points in $\bx^\flat \setminus \by^\flat$, and
denote by $n_{1,1}$, $\ldots$, $n_{1,m_1}$ the weights prescribed
to the points of $\bx^\flat \setminus \by^\flat$ by the sequence
$\alp$. Furthermore, put $m_2 = \|\bet' - \bet\|$, consider a
non-increasing sequence $y_{2,1}$, $\ldots$, $y_{2,m_2}$ of second
coordinates of the images of terminal univalent vertices of left
ends belonging to~${\cal E}$, and denote by $n_{2,1}$, $\ldots$,
$n_{2,m_2}$ the weights of the corresponding left ends. Take a
parameterized plane tropical curve $(\Gam,w,h)$ satisfying the
conditions~(1)-(4) of Lemma~\ref{lch26}. Due to a possibility to
make a negative horizontal shift of~$p$ and simultaneously
compose~$h$ with this shift (see the definition of~${\cal X}$), we
can assume that there exists a vertical line $x = c$ such that the
images of vertices of~$\Gam$ lie in the left half-plane delimited
by this line, and the images of non-univalent vertices of
$\coprod_\newi\Gam^{(\newi)}$ lie in the right half-plane. In
particular, the line $x=c$ crosses the images of all left ends
belonging to~${\cal E}$ and the horizontal positively directed
images of ends of~$\Gam$.

Cut along the preimages of $x = c$ the left ends belonging to
${\cal E}$ and the ends of~$\Gam$ whose images are horizontal and
positively directed, and remove the trivial pieces of the edges
cut. The natural gluing of remaining pieces of $\Gam$ and
$\coprod_{\newi=1}^\newm\Gam^{(\newi)}$ gives rise to a collection
of parameterized plane tropical curves $(\widetilde\Gam^{(\newj)},
\widetilde{w}^{(\newj)}, \widetilde{h}^{(\newj)})$, $\newj = 1$,
$\ldots$, $\newl$, a splitting
$(\Del,\alp,\bet,g)=\sum_{\newj=1}^\newl
(\widetilde\Del^{(\newj)},
\widetilde\alp^{(\newj)},\widetilde\bet^{(\newj)},
\widetilde{g}^{(\newj)})$ in ${\cal S}$, and a partition
$\bigsqcup_{\newj=1}^\newl\widetilde\bx^{(\newj)}$ of~$\bx$ which
satisfy the following properties:
\begin{itemize}
\item each curve $(\widetilde\Gam^{(\newj)}, \widetilde{w}^{(\newj)},
\widetilde{h}^{(\newj)})$ matches
the configuration $\widetilde\bx^{(\newj)}$ and has
the degree
$(\widetilde\Del^{(\newj)})^{\widetilde\alp^{(\newj)},
\widetilde\bet^{(\newj)}}$
and the genus~$\widetilde{g}^{(\newj)}$,
\item the partition $\bigsqcup_{\newj=1}^\newl\widetilde\bx^{(\newj)}$
is compatible with the splitting
$(\Del,\alp,\bet,g)=\sum_{\newj=1}^\newl
(\widetilde\Del^{(\newj)},\widetilde\alp^{(\newj)},
\widetilde\bet^{(\newj)},
\widetilde{g}^{(\newj)})$.
\end{itemize}
Since the configuration~$\bx$ is multi-tropically generic, each
configuration~$\widetilde\bx^{(j)}$, $\newj = 1$, $\ldots$,
$\newl$ is tropically generic. This implies, that for any $\newj =
1$, $\ldots$, $\newl$ the curve $(\widetilde\Gam^{(\newj)},
\widetilde{w}^{(\newj)}, \widetilde{h}^{(\newj)})$ defines a nodal
plane tropical curve~$\widetilde{T}^{(\newj)}$ such that no point
in $(\widetilde\bx^{(\newj)})^\sharp$ is a vertex
of~$\widetilde{T}^{(\newj)}$. Moreover, according to
Lemma~\ref{good-points}, the lift~$\widetilde\bpp^{(\newj)}$
of~$\widetilde\bx^{(\newj)}$ to $\widetilde\Gam^{(\newj)}$
produces a parameterized marked tropical curve
$(\widetilde\Gam^{(\newj)}, \widetilde{w}^{(\newj)},
\widetilde{h}^{(\newj)}, \widetilde\bpp^{(\newj)})$, and thus, the
sum $\widetilde{T}^{(1)} + \ldots + \widetilde{T}^{(\newl)}$ is a
nodal plane tropical curve. Hence, the configuration~$\by$ is
multi-tropically generic. \proofend

\subsection{Proof of Theorem~\ref{formula}}
Choose positive real numbers~$\varepsilon$ and~$N$ satisfying the
inequality $\varepsilon < N$, and consider a configuration $\bx =
(\bx^\flat, \bx^\sharp) \in {\cal X} \subset
\Omega(\Del,\alp,\bet,g)$ such that the second coordinates of all
the points in~$\bx$ belong to the interval $(-\varepsilon,
\varepsilon)$, the first coordinate of one point $p \in
\bx^\sharp$ is smaller than $-N$, and the first coordinates of the
other points in~$\bx^\sharp$ belong to the interval
$(-\varepsilon, \varepsilon)$. Assume that the
numbers~$\varepsilon$ and~$N$ are chosen in such a way that, for
some numbers~$a$ and~$b$ satisfying the inequalities $-N < a < b <
-\varepsilon$, the intersection of any tropical curve in ${\cal
T}(\Del,\alpha, \beta, g, \bx)$ with the rectangle $R_a^b = \{(x,
y) \in \R^2: \; a \leq x \leq b \quad \text{\rm and} \;
-\varepsilon \leq y \leq \varepsilon\}$ consists of horizontal
segments (the existence of such~$\varepsilon$ and~$N$ is
guaranteed by Lemma~\ref{splitting}). Consider a tropical curve $T
\in {\cal T}(\Del,\alpha, \beta, g, \bx)$ without edges of even
weight.

Suppose that~$p$ belongs to a $\sig$-end~$e$ of~$T$. Since the
configuration $\bx$ is multi-tropically generic,
Lemma~\ref{one-lemma} implies that the end~$e$ is not marked, {\it
i.e.}, terminates at a point~$q \in L_{-\infty}$ different from
any point of~$\bx^\flat$. Consider the configuration $\widehat\bx
= (\widehat\bx^\flat, \widehat\bx^\sharp) \in \Omega(\Del, \alp +
\theta_k, \bet - \theta_k, g)$, where~$k=w(e)$ is the weight
of~$e$, the configuration~$\widehat\bx^\flat$ is obtained
from~$\bx^\flat$ by insertion of~$q$ in the $k$-th group of
points, and $\widehat\bx^\sharp = \bx^\sharp \setminus \{p\}$.
Since $\bx \in {\cal X}$, the configuration~$\widehat\bx$ is
multi-tropically generic.

Clearly, the curve~$T$ belongs to ${\cal T}(\Del,\alpha +
\theta_k, \beta - \theta_k, g, \widehat\bx)$. On the other hand,
we can assume that~$N$ is chosen so that ${\cal T}(\Del,\alpha +
\theta_k, \beta - \theta_k, g, \widehat\bx) \subset {\cal
T}(\Del,\alpha, \beta, g, \bx)$. Thus, the contribution to
$W(\Del,\alp,\bet,g)$ of the curves $T \in {\cal T}(\Del,\alpha,
\beta, g, \bx)$ such that~$p$ belongs to a $\sig$-end of~$T$ is
equal to
$$\sum_{\renewcommand{\arraystretch}{0.6}
\begin{array}{c}
\scriptstyle{k\ge 1}\\
\scriptstyle{\beta_k>0}
\end{array}}
W(\Del,\alpha+\theta_k, \beta-\theta_k, g).$$

Suppose now that~$p$ does not belong to any $\sig$-end of~$T$.
Since~$T$ is nodal, it can be represented in a unique way as a sum
of its irreducible subcurves $T^{(1)}$, $\ldots$, $T^{(\newl)}$.
For any $\newj = 1$, $\ldots$, $\newl$, consider a simple
parameterization $(\Gam^{(\newj)}, w^{(\newj)}, h^{(\newj)})$
of~$T^{(\newj)}$. Pick a number~$c$ such that $a < c < b$. For any
curve $(\Gam^{(\newj)}, w^{(\newj)}, h^{(\newj)})$ consider the
lift~$\Upsilon^{(\newj)} \subset \Gam^{(\newj)}$ of the
intersection points of~$T^{(\newj)}$ with the vertical segment $x
= c$, $-\varepsilon \leq y \leq \varepsilon$.
Lemma~\ref{splitting} implies that no connected component of
$\Gam^{(\newj)} \setminus \Upsilon^{(j)}$ has an image
intersecting the both halves $R_a^b \cap \{x < c\}$ and $R_a^b
\cap \{x > c\}$ of~$R_a^b$. If $\Upsilon^{(\newj)}\ne
\varnothing$, then $\Upsilon^{(\newj)}$ cuts $\Gam^{(\newj)}$ in
two parts: the image of any connected component in the right part
$\Gam^{(\newj)}_{\cal R}$ intersects $R_a^b \cap \{x > c\}$, and
the image of any connected component of the left part
$\Gam^{(\newj)}_{\cal L}$ intersects $R_a^b \cap \{x < c\}$. Any
connected component~$\Gam^{(\newj, \nu)}_{\cal L}$, $\nu = 1$,
$\ldots$, $\newl_j$, gives rise to a parameterized plane tropical
curve $(\Gam^{(\newj,\nu)}_{\cal L}, w^{(\newj,\nu)}_{\cal L},
h^{(\newj,\nu)}_{\cal L})$, where the weight function $w^{(\newj,
\nu)}_{\cal L}$ is induced by $w^{(\newj)}$, and $h^{(\newj,
\nu)}_{\cal L}$ is given by a modification of~$h^{(\newj)}$ on the
edges cut (without changing the directions of the images of these
edges). The image of~$\Gam^{(\newj,\nu)}_{\cal L}$ under
$h^{(\newj,\nu)}_{\cal L}$ is obtained from
$h^{(\newj)}(\Gam^{(\newj,\nu)}_{\cal L})$ by the extension of the
edges cut by the segment $x = c$, $-\varepsilon \leq y \leq
\varepsilon$ up to horizontal positively directed rays.

For any connected component $\Gam^{(\newj,\mu)}_{\cal R}$, $\mu =
1$, $\ldots$, $\newl'_\newj$ of $\Gam^{(\newj)}_{\cal R}$, denote
by $\widehat\Gam^{(\newj,\mu)}_{\cal R}$ the graph obtained from
$\Gam^{(\newj,\mu)}_{\cal R}$ by adding a vertex to each cut edge.
The graph $\widehat\Gam^{(\newj,\mu)}_{\cal R}$ gives rise to a
parameterized plane tropical curve
$(\widehat\Gam^{(\newj,\mu)}_{\cal R}, w^{(\newj,\mu)}_{\cal R},
h^{(\newj,\mu)}_{\cal R})$, where again the weight function
$w^{(\newj,\mu)}_{\cal R}$ is induced by $w^{(\newj)}$, the
restriction of $h^{(\newj,\mu)}_{\cal R}$ on
$\Gam^{(\newj,\mu)}_{\cal R}$ is given by a modification
of~$h^{(\newj)}$ on the edges cut (without changing the directions
of the images of these edges), and the image under
$h^{(\newj,\mu)}_{\cal R}$ of any added vertex belongs to
$L_{-\infty}$. The image of~$\widehat\Gam^{(\newj,\mu)}_{\cal R}$
under $h^{(\newj,\mu)}_{\cal R}$ is obtained from
$h^{(\newj)}(\Gam^{(\newj,\mu)}_{\cal R})$ by the extension of the
edges cut by the segment $x = c$, $-\varepsilon \leq y \leq
\varepsilon$ up to horizontal negatively directed rays. Denote by
$T'$ the plane tropical curve defined by the collection of all
parameterized plane tropical curves
$(\widehat\Gam^{(\newj,\mu)}_{\cal R}, w^{(\newj,\mu)}_{\cal R},
h^{(\newj,\mu)}_{\cal R})$, $\newj = 1$, $\ldots$, $\newl$, $\mu =
1$, $\ldots$, $\newl'_\newj$. We say that~$T'$ is the {\it
derivation} of~$T$.

If $\Upsilon^{(\newj)}$ is empty,
then $h^{(\newj)}(\Gam^{(\newj)}) \subset \{x < c\}$,
and we use
the notation
$(\Gam^{(\newj,1)}_{\cal L}, w^{(\newj,1)}_{\cal L}, h^{(\newj,1)}_{\cal L})$
for the parameterized plane tropical curve
$(\Gam^{(\newj)}, w^{(\newj)}, h^{(\newj)})$.

Since the half-plane $x < c$ contains only one point of
$\bx^\sharp$, Lemma~\ref{good-points} implies that
\begin{itemize}
\item all the curves
$(\Gam^{(\newj,\nu)}_{\cal L},
w^{(\newj,\nu)}_{\cal L}, h^{(\newj,\nu)}_{\cal L})$,
$\newj = 1$, $\ldots$, $\newl$, $\nu = 1$, $\ldots$, $\newl_\newj$,
are of genus~$0$,
\item among the curves
$(\Gam^{(\newj,\nu)}_{\cal L}, w^{(\newj,\nu)}_{\cal L},
h^{(\newj,\nu)}_{\cal L})$, $\newj = 1$, $\ldots$, $\newl$, $\nu =
1$, $\ldots$, $\newl_\newj$, there exists a curve $(\Gam, w, h)$
such that~$\Gam$ has exactly one end whose image under~$h$ points
upwards ({\it i.e.}, has the direction either $(0, 1)$ or $(1,
1)$), and exactly one end whose image points downwards ({\it
i.e.}, has the direction either $(0, -1)$ or $(-1, -1)$),
\item the image under~$h$ of any
left end of~$\Gam$ terminates at a point of~$\bx^\flat$,
\item for any curve $(\Gam^{(\newj,\nu)}_{\cal L},
w^{(\newj,\nu)}_{\cal L}, h^{(\newj,\nu)}_{\cal L})$ different from
$(\Gam, w, h)$ the image $h^{(\newj,\nu)}_{\cal
L}(\Gam^{(\newj,\nu)}_{\cal L})$ is a horizontal straight line; such
curves $(\Gam^{(\newj,\nu)}_{\cal L}, w^{(\newj,\nu)}_{\cal L},
h^{(\newj,\nu)}_{\cal L})$ are called {\it horizontal}.
\end{itemize}
In particular,
the Newton polygon~$\Del'$ of~$T'$ is the
$\daleth$-peeling $l_\daleth(\Del)$
of~$\Del$ for
certain $\daleth \subset \{\overrightarrow{0}, \overrightarrow{-1}\}$.
Denote by $\sig'$ the vertical left-most side
of~$\Del'$.

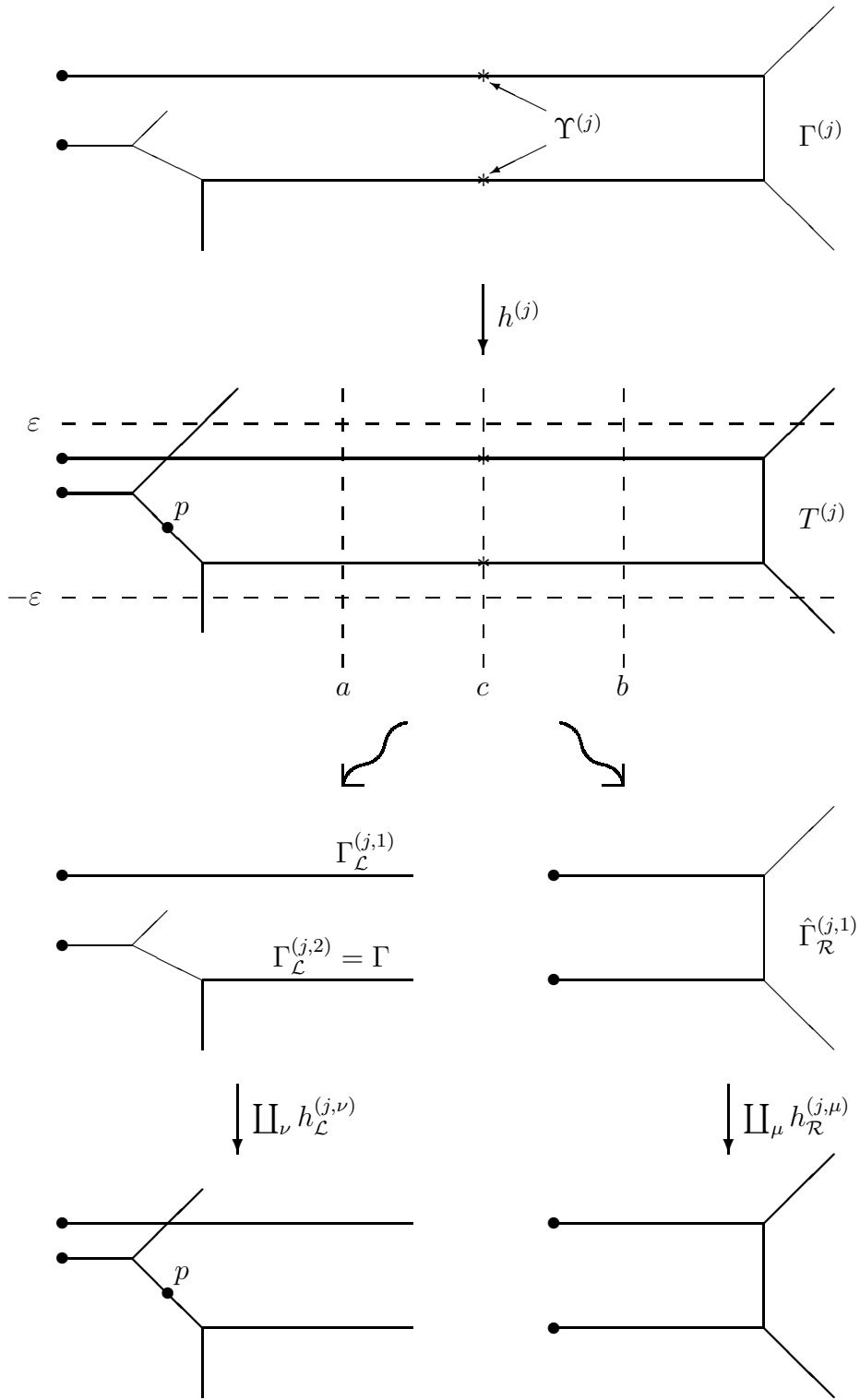
\begin{figure}
\setlength{\unitlength}{1cm}
\begin{picture}(12,20)(-1,0)
\thinlines\put(0,19){\line(1,0){10}}\put(10,19){\line(1,1){1}}
\put(10,19){\line(0,-1){1.5}}\put(10,17.5){\line(1,-1){1}}
\put(10,17.5){\line(-1,0){8}}\put(2,17.5){\line(0,-1){1}}
\put(2,17.5){\line(-2,1){1}}\put(1,18){\line(-1,0){1}}
\put(1,18){\line(1,1){0.5}}\thicklines\put(0,13.5){\line(1,0){10}}
\put(10,13.5){\line(1,1){1}}\put(10,13.5){\line(0,-1){1.5}}
\put(10,12){\line(1,-1){1}}\put(10,12){\line(-1,0){8}}
\put(2,12){\line(0,-1){1}}\put(2,12){\line(-1,1){1}}
\put(1,13){\line(-1,0){1}}\put(1,13){\line(1,1){1.5}}\thinlines
\put(0,7.5){\line(1,0){5}}\put(1,6.5){\line(2,-1){1}}\put(1,6.5){\line(1,1){0.5}}
\put(0,6.5){\line(1,0){1}}\put(2,6){\line(0,-1){1}}
\put(2,6){\line(1,0){3}}\put(7,7.5){\line(1,0){3}}
\put(7,6){\line(1,0){3}}\put(10,6){\line(0,1){1.5}}
\put(10,6){\line(1,-1){1}}\put(10,7.5){\line(1,1){1}}\thicklines
\put(0,2.5){\line(1,0){5}}\put(0,2){\line(1,0){1}}
\put(1,2){\line(1,-1){1}}\put(2,1){\line(0,-1){1}}
\put(2,1){\line(1,0){3}}\put(1,2){\line(1,1){1}}
\put(7,2.5){\line(1,0){3}}\put(7,1){\line(1,0){3}}
\put(10,1){\line(0,1){1.5}}\put(10,1){\line(1,-1){1}}
\put(10,2.5){\line(1,1){1}}\thinlines \dashline{0.2}(0,14)(11,14)
\dashline{0.2}(0,11.5)(11,11.5)\dashline{0.2}(4,10.5)(4,14.5)
\dashline{0.2}(6,10.5)(6,14.5)\dashline{0.2}(8,10.5)(8,14.5)
\put(-0.1,18.9){$\bullet$}\put(-0.1,17.9){$\bullet$}
\put(-0.1,13.4){$\bullet$}\put(-0.1,12.9){$\bullet$}
\put(-0.1,7.4){$\bullet$}\put(-0.1,6.4){$\bullet$}
\put(-0.1,2.4){$\bullet$}\put(-0.1,1.9){$\bullet$}
\put(6.9,5.9){$\bullet$}\put(6.9,7.4){$\bullet$}
\put(6.9,2.4){$\bullet$}\put(6.9,0.9){$\bullet$}
\put(1.4,12.4){$\bullet$}\put(1.4,1.4){$\bullet$}
\put(5.9,11.9){$*$}\put(5.9,13.4){$*$}\put(5.9,17.4){$*$}
\put(5.9,18.9){$*$}\put(10.5,18){$\Gamma^{(j)}$}\put(10.5,12.5){$T^{(j)}$}
\put(-0.5,13.9){$\eps$}\put(-0.8,11.4){$-\eps$}\put(3.9,10.1){$a$}
\put(5.9,10.1){$c$}\put(7.9,10.1){$b$}\put(1.6,12.7){$p$}
\put(1.6,1.7){$p$}\put(3.9,7.7){$\Gam^{(j,1)}_{\cal L}$}
\put(3,6.2){$\Gam^{(j,2)}_{\cal
L}=\Gam$}\put(10.5,6.5){$\hat\Gam^{(j,1)}_{\cal
R}$}\put(7,18.1){$\Upsilon^{(j)}$}\put(6.9,18.5){\vector(-2,1){0.8}}
\put(6.9,18){\vector(-2,-1){0.8}}
\thicklines\put(6,16){\vector(0,-1){1}}\put(6.2,15.4){$h^{(j)}$}
\put(2.5,4.5){\vector(0,-1){1}}\put(2.7,3.9){$\coprod_\nu
h^{(j,\nu)}_{\cal
L}$}\put(9.5,4.5){\vector(0,-1){1}}\put(9.7,3.9){$\coprod_\mu
h^{(j,\mu)}_{\cal R}$} \thicklines
\put(4,8.8){\line(1,0){0.3}}\put(4,8.8){\line(0,1){0.3}}
\qbezier(4,8.8)(4.05,9.05)(4.3,9.1)\qbezier(4.3,9.1)(4.55,9.15)(4.6,9.4)\qbezier(4.6,9.4)(4.65,9.65)(4.9,9.7)
\qbezier(8,8.8)(7.95,9.05)(7.7,9.1)\qbezier(7.7,9.1)(7.45,9.15)(7.4,9.4)\qbezier(7.4,9.4)(7.35,9.65)(7.1,9.7)
\put(8,8.8){\line(-1,0){0.3}}\put(8,8.8){\line(0,1){0.3}}
\end{picture}
\caption{Cut in the proof of Theorem~\ref{formula}}\label{chfig8}
\end{figure}

Let $(\bx')^\flat \subset \bx^\flat$ be the subconfiguration
formed by the images of marked terminal points of horizontal
curves, and let $\alpha' \leq \alp$ be the corresponding sequence
in ${\cal C}$. Since~$T$ is nodal, the images of terminal points
of left ends of graphs $\Gam^{(\newj)}_{\cal R}$ are disjoint from
the images of terminal points of left ends of~$\Gam$. Among the
left ends of graphs $\Gam^{(\newj)}_{\cal R}$ consider the edges
whose images terminate at points different from the points
of~$\bx^\flat$, and denote by $\beta'$ the sequence determined by
the weights of the edges considered. Since the image under~$h$ of
any left end of~$\Gam$ terminates at a point of~$\bx^\flat$, we
have $\beta' \geq \beta$. Furthermore, $J\alpha' + J\beta' =
|\sig'|$. Counting the edges cut by the segment $x = c$,
$-\varepsilon \leq y \leq \varepsilon$, we obtain that the
genus~$g'$ of~$T'$ is equal to $g - ||\beta' - \beta|| + 1$. The
curve~$T'$ belongs to ${\cal T}(\Del',\alpha',\beta', g', \bx')$,
where $\bx' = ((\bx')^\flat, (\bx')^\sharp)$, and $(\bx')^\sharp$
is obtained from $\bx^\sharp$ by removing the point~$p$, and all
the edges of~$T'$ are of odd weights.

Describe now the inverse procedure. Fix two elements $\alpha',
\beta' \in {\cal C}$ such that $\alpha' \leq \alpha$, $\beta' \geq
\beta$, and $J\alpha' + J\beta' = |\sig'|$. Put $g' = g - ||\beta'
- \beta|| + 1$. Choose a subconfiguration $(\bx')^\flat \subset
\bx^\flat$ such that $(\bx')^\flat$ corresponds to the sequence
$\alpha'$ (the number of such subconfigurations $(\bx')^\flat$ is
$\left(
\begin{matrix}\alpha\\ \alpha'\end{matrix}\right)$).
Consider the configuration $\bx' = (\bx')^\flat \cup (\bx')^\sharp$,
where $(\bx')^\sharp$, as before, is
obtained from $\bx^\sharp$ by removing the point~$p$.
By Lemma
\ref{lch27},
the configuration~$\bx'$ is multi-tropically generic.
Furthermore, by Lemma~\ref{multi-finiteness} the
set ${\cal T}(\Del', \alpha', \beta', g', \bx')$ is finite,
and we
can assume that the first coordinate of~$p$ is much less than the
first coordinates of the vertices of all the curves
in ${\cal T}(\Del',\alpha', \beta', g', \bx')$.

Pick
a curve $T' \in {\cal T}(\Del',\alpha', \beta', g', \bx')$
without edges of even weight.
Let $(\widehat\Gam^{(\newi)}_{\cal R},
w^{(\newi)}_{\cal R}, h^{(\newi)}_{\cal R})$,
$\newi = 1$, $\ldots$, $\newm$, be simple parameterizations
of irreducible subcurves of~$T'$.
Among the non-marked left ends
of $\coprod_{\newi = 1}^{\newm}\widehat\Gam^{(\newi)}_{\cal R}$
choose $||\bet'-\bet||$ ends
whose weights fit the sequence
$\bet'-\bet$
(the
number of such choices
is
$\left(
\begin{matrix}\beta'\\ \beta\end{matrix}\right)$),
and denote the chosen set by~${\cal E}$. Lemma~\ref{lch26}
provides a parameterized plane tropical curve $(\Gam, w, h)$ such
that
\begin{itemize}
\item $(\Gam, w, h )$ is of genus~$0$ and
has exactly two ends whose images under~$h$ are not horizontal;
these ends are of weight~$1$, and their images
point in the directions of outward normal vectors
of
the sides $\top(\Del, \daleth)$ and $\bot(\Del, \daleth)$
of~$\Del$
(see the proof of Lemma~\ref{lch27} for notation),
\item the images of left ends of~$\Gam$
terminate at the points of $\bx^\flat\backslash(\bx')^\flat$,
and the weights of these left ends are given
by the sequence $\alp - \alp'$,
\item the ends of~$\Gam$ whose images
are horizontal positively directed fit the
ends belonging to~${\cal E}$ and have
the corresponding weights,
\item $p\in h(\Gam)$.
\end{itemize}
In the same way as in the proof of Lemma \ref{lch27}, we can glue
$(\Gam, w, h)$ with the curves $(\widehat\Gam^{(\newi)}_{\cal R},
w^{(\newi)}_{\cal R}, h^{(\newi)}_{\cal R})$, $\newi = 1$,
$\ldots$, $\newm$, and obtain a collection of parameterizations of
plane tropical curves $T^{(1)}$, $\ldots$, $T^{(\newl)}$ whose
sum~$T$ belongs to~${\cal T}(\Del, \alp, \bet, g, \bx)$.
Since~$\bx$ is multi-tropically generic, Lemma~\ref{lch26} implies
that the curve~$(\Gam, w, h)$ is defined by the above properties
uniquely. Furthermore, the initial curve~$T'$ is the derivation
of~$T$. The multiplicative contribution of the trivalent vertices
of~$\Gam$ to the Welschinger multiplicity of~$T$ is~$1$, and we
finally conclude that, for given $\alpha'$ and $\beta'$, the
contribution to $W(\Del,\alp,\bet,g)$ of the curves $T \in {\cal
T}(\Del,\alpha, \beta, g, \bx)$ such that~$p$ does not belong to
any $\sig$-end of~$T$ is equal to $ \left(
\begin{matrix}\alpha\\ \alpha'\end{matrix}\right)\left(\begin{matrix}\beta'\\
\beta\end{matrix} \right)W(\Del',\alpha', \beta', g')$. \proofend

\subsection{Proof of Theorem \ref{cht3}}
The proof goes in the same way as for Theorem~\ref{formula}. The
only modification concerns the second sum in the right-hand side
of (\ref{che33}): assuming that a plane tropical curve $T'\in
{\cal T}(\Del', \alp', \bet', g', \bx')$ is the derivation of
$T\in{\cal T}^{\irr}(\Del, \alp, \bet, g, \bx)$, we have to
describe possible irreducible subcurves of~$T'$ and their
contribution to the formula. The first (respectively, second)
coefficient in the second sum of the formula reflects the
distribution of the points of $(\bx')^\flat$ (respectively,
$(\bx')^\sharp$) among the irreducible subcurves of~$T'$. The
conditions on the
numbers $\alp',\bet',g'$ and $\alp^{(i)},\bet^{(i)},g^{(i)}$,
$i=1,...,m$, come from the conditions in Theorem~\ref{formula},
and the inequalities $||\widetilde\bet^{(i)}||>0$ mean that, for
each irreducible subcurve~$T^{(i)}$, its simple parameterization
$(\widehat\Gam^{(i)}_{\cal R}, w^{(i)}_{\cal R}, h^{(i)}_{\cal
R})$ must glue with the curve $(\Gam, w, h)$ (in the notation of
the proof of Theorem~\ref{formula}). \proofend

\section{Concluding remarks}\label{remarks}

\subsection{Generating functions}\label{chsec5}

Let~$\Del$ be one of the polygons shown
in Figure~\ref{chf4},
and~$\Sig$ the real toric Del Pezzo surface defined by~$\Del$. We
say that a convex lattice polygon $\Del^0$ has the same {\it
shape} as $\Del$, if $\Del$ and $\Del^0$ have the same number of
sides, and any side of~$\Del^0$ is parallel to a side of~$\Del$.
Let $\Xi_\Sig$ be the set which consists of all
convex lattice polygons (considered up to parallel translation)
having the same shape as~$\Del$ and their
$\daleth$-peelings
($\daleth \subset \{\overrightarrow{0}, \overrightarrow{-1}\}$
being
$\Del$-admissible).
The set $\Xi_\Sig$ is a commutative
semigroup with respect to the Minkowsky sum.

Following \cite{Ge,Va} introduce two generating functions
$$Z_\Sig(w,x,y,z)=\sum_{\renewcommand{\arraystretch}{0.6}
\begin{array}{c}
\scriptstyle{g\in\Z,\ \Del\in\Xi_\Sig}\\
\scriptstyle{\alp,\bet\in{\cal C},\ J\alp + J\bet = |\sig(\Del)|}
\end{array}}W(\Del,\alp,\bet,g)v^\Del w^{g-1} \ \frac{x^\alp}{\alp!}
\ y^\bet \ \frac{z^r}{r!} \ ,$$
$$Z^{\irr}_\Sig(w,x,y,z)=\sum_{\renewcommand{\arraystretch}{0.6}
\begin{array}{c}
\scriptstyle{g\ge 0,\ \Del\in\Xi_\Sig}\\
\scriptstyle{\alp,\bet\in{\cal C},\ J\alp + J\bet = |\sig(\Del)|}
\end{array}}W^{\irr}(\Del,\alp,\bet,g)
v^\Del w^{g-1} \ \frac{x^\alp}{\alp!} \ y^\bet \ \frac{z^r}{r!} \
,$$ where~$\sig(\Del)$ is the intersection of~$\Del$ with its left
vertical supporting line, $x$ and~$y$ are infinite sequences of
variables, and~$r$ is defined by~(\ref{che201}). These generating
functions can be seen as formal series in variables~$w,z$ and
multi-variables $x,y$ with coefficients in the Novikov ring of the
semigroup~$\Xi_\Sig$. Using the same arguments as in \cite{Va},
Section 6.4 and \cite{Ge}, Section 5.3, one can check that these
generating functions are related by the identity $Z_\Sig=\exp
Z^{\irr}_\Sig$ and satisfy the differential equations

$$\sum_{\daleth \subset \{\overrightarrow{0}, \overrightarrow{-1}\}}v^{l_\daleth}\left(\frac{\partial}
{\partial z}-\sum_{k=1}^\infty y_k\frac{\partial}{\partial
x_k}\right)Z_\Sig=\frac{1}{w}\cdot
\res_{t=0}\exp\sum_{k=1}^\infty\left(t^{-k}x_k+wt^k\frac{\partial}{\partial
y_k}\right)Z_\Sig\ ,$$
$$\sum_{\daleth \subset \{\overrightarrow{0}, \overrightarrow{-1}\}}v^{l_\daleth}\left(\frac{\partial}
{\partial z}-\sum_{k=1}^\infty y_k\frac{\partial}{\partial
x_k}\right)Z^{\irr}_\Sig$$
$$=\frac{1}{w}\cdot\res_{t=0}\exp\left(\sum_{k=1}^\infty
(t^{-k}x_k+Z^{\irr}_\Sig\big|_{y_k\mapsto
y_k+wt^k})-Z^{\irr}_\Sig\right)\ ,$$ where
$Z^{\irr}_\Sig\big|_{y_k\mapsto y_k+wt^k}$ stands for
$Z^{\irr}_\Sig$ with $y_k$ replaced by $y_k+wt^k$, and
$$
v^{l_{\daleth}}v^\Del= \begin{cases} v^{l_{\daleth}(\Del)}, \quad &
\text{if} \; \daleth \; \text{is}\: \Del-\text{admissible},
 \\
0, \quad & \text{otherwise}.
\end{cases}
$$

\subsection{Non-invariance in the classical
setting}\label{noninv}
The absence of
invariants of topological nature as mentioned in Introduction can
be illustrated by the following examples.

Real irreducible plane curves of degree~$d$ and genus~$g$
passing through a generic configuration of $3d+g-1$ real points
form a finite set.
If $g > 0$, then under variation of the point configuration,
this set is subject to codimension one events in which a pair of
real curves with the same embedded topology disappears turning
into a pair of imaginary conjugate curves, or vice versa. In the
case of elliptic quartics ($d = 4$, $g = 1$),
it is shown in \cite{IKS}, Theorem 3.1.
Gluing elliptic quartics with one or several lines, it is not
difficult to construct examples of higher degree and genus.
The example with elliptic quartics shows also
that the number $W(\Del(5\PP^1),(0),(5),0)$ does
not lift up to the classical setting
as an invariant formulated in purely
topological terms.

Another example, demonstrating the same phenomenon, is as follows.
Consider plane rational curves of degree $d$ which pass through
$3d+1-b_1-b_2$ generic points outside a line $\PP^1 \subset \PP^2$
and have two non-fixed tangency points with $\PP^1$, one of
intersection order~$b_1$ and the other of intersection
order~$b_2$, where $b_1$ and $b_2$ are distinct, odd, and satisfy
the inequality $b_1+b_2<d$. Under variation of the point
configuration, the collision of the two tangency points into one
tangency point of order $b_1 + b_2$ is an event of codimension
one. Crossing such a wall leads to appearance, or disappearance,
of two real curves which have the same embedded topology (even
with respect to $\PP^1$).

One more phenomenon is the change of the Welschinger multiplicity
of precisely one member of the set of curves. Consider real plane
rational curves of degree $d\ge 5$ passing through $a \geq 3$
fixed generic real points on $\PP^1 \subset \PP^2$ and $3d - 1 -
a$ generic real points in $\PP^2 \setminus \PP^1$. Here we observe
the following codimension one event: precisely one of the curves
splits into $D+\PP^1$, where $D$ is tangent to $\PP^1$ at one
point and transversal to it at $d-3$ other points. On one side of
such a wall the tangency point turns into a solitary node, and on
the other side into a crossing point, whereas the other
singularities do not change.

The above arguments do not exclude the existence of relative real algebraic
enumerative invariants in other situations. For example, such invariants
were introduced by Welschinger \cite{W3} in the case of one simple tangency
constraint with respect to a smooth null-homologous curve.

{\ncsc Universit\'{e} Louis Pasteur et IRMA \\[-21pt]

7, rue Ren\'{e} Descartes, 67084 Strasbourg Cedex, France} \\[-21pt]

{\it E-mail}: {\ntt itenberg@math.u-strasbg.fr}

\vskip10pt

{\ncsc Universit\'{e} Louis Pasteur et IRMA \\[-21pt]

7, rue Ren\'{e} Descartes, 67084 Strasbourg Cedex, France} \\[-21pt]

{\it E-mail}: {\ntt kharlam@math.u-strasbg.fr}

\vskip10pt

{\ncsc School of Mathematical Sciences \\[-21pt]

Raymond and Beverly Sackler Faculty of Exact Sciences\\[-21pt]

Tel Aviv University,
Ramat Aviv, 69978 Tel Aviv, Israel} \\[-21pt]

{\it E-mail}: {\ntt shustin@post.tau.ac.il}

\end{document}